\documentclass[10pt]{siamltex}
\usepackage{stmaryrd}
\usepackage{mathrsfs}%
\usepackage{amsfonts}
\usepackage{amsmath}
\usepackage{amssymb}
\usepackage{graphicx}
\usepackage{latexsym}
\usepackage{indentfirst}
\usepackage{graphicx}
\usepackage{overpic}
\usepackage{enumerate}
\usepackage{slashbox}
\usepackage[bf]{caption2}
\usepackage{cite}
\usepackage{bbding}
\usepackage{multirow}%

\newtheorem{thm}{Theorem}
\newtheorem{remark}{Remark}[section]
\newtheorem{lem}{Lemma}

\renewcommand{\div}{\operatorname{div}}

\begin{document}

\title{An $hp$-version error analysis of the discontinuous Galerkin method for linear elasticity}%
\author{
Jianguo Huang\thanks{School of Mathematical Science, and MOE-LSC, Shanghai Jiao Tong University, Shanghai 200240, China (jghuang@sjtu.edu.cn). The work of the first author was partly supported by NSFC (Grant nos. 11571237 and 11171219).}
\and
Xuehai Huang\thanks{Corresponding author. College of Mathematics and Information Science, Wenzhou University, Wenzhou 325035, China (xuehaihuang@wzu.edu.cn). The work of this author was partly supported by NSFC (Grant nos. 11771338 and 11301396), Zhejiang Provincial Natural Science Foundation of China Project (Grant no. LY17A010010), and Wenzhou Science and Technology Plan Project (Grant no. G20160019).}
}
%
%

\maketitle

\noindent\textbf{Abstract:}
An $hp$-version error analysis is developed for the general DG method in mixed formulation for solving the
linear elastic problem. First of all, we give the $hp$-version error estimates of two $L^2$ projection operators.
Then incorporated with the techniques in \cite{CastilloCockburnPerugiaSchotzau2000}, we obtain the $hp$-version error estimates in energy norm and $L^2$ norm. Some numerical experiments are provided for demonstrating the theoretical results.

\vspace{1 pc}

\vspace{1 pc}

\vspace{1 pc}

\section{Introduction}

The linear elastic equations are used to describe the deformation of elastic structures under the action of prescribed loads, which are the fundamental equations in mathematical physics. Both the displacement and stress fields are the fundamental physical quantities in mechanical analysis. When approximating the displacement and stress simultaneously, the mixed finite element method (cf. \cite{ArnoldWinther2002,AdamsCockburn2005, ArnoldDouglasGupta1984, HuManZhang2014, Hu2015a,HuZhang2015,HuZhang2016, ChenHuHuang2017}) can achieve higher-accuracy stress than the standard displacement finite element method, with which the stress is obtained by differentiating the displacement and using the constitutive law of stress-strain. The main and critical difficulty in construction of such mixed methods is largely due to the fact that the stress tensor is  symmetric and belong to $\boldsymbol{H}(\textrm{div}):=\{\boldsymbol{\tau}\in (L^2(\Omega))_{d\times d}; \mathbf{div}~\boldsymbol{\tau}\in (L^2(\Omega))^d\}$.
To overcome his difficulty, the discontinuous Galerkin (DG) finite element method is an apt choice to solve the linear elasticity by weakening the regularity of finite element space. Historically, some discontinuous Galerkin methods for linear elasticity are presented in primal formulation, including the local DG (LDG) method in \cite{LewNeffSulskyOrtiz2004}, the compact DG method in \cite{HuangHuang2013a}, and the interior penalty DG method in \cite{HansboLarson2002,HansboLarson2003}. A mixed discontinuous Galerkin finite element method with symmetric stress tensor for linear elasticity is given in \cite{CaiYe2005}, which is a special case of general DG formulation in \cite{ChenHuangHuangXu2010}.
More recently, some stabilized mixed finite element methods with symmetric stress tensor for linear elasticity, based on the Hu-Zhang element in \cite{Hu2015a,HuZhang2015c,HuZhang2015}, are introduced in \cite{ChenHuHuang2017}.
In \cite{Bustinza2006}, a mixed DG formulation is also designed but the stress tensor is nonsymmetric.
Following the ideas in \cite{CastilloCockburnPerugiaSchotzau2000, Cockburn2003}, a general framework of constructing
DG methods with symmetric stress tensor has been developed in \cite{ChenHuangHuangXu2010}   for solving the linear elasticity problem, and the $h$-optimal convergence of the resulting LDG method is developed as well.

Polynomials of arbitrary degree can be taken on each element in the discontinuous Galerkin method, for the continuity of finite element spaces across the interfaces of triangulation is not required. Thus it is natural to analyze DG method in $hp$-version context. To this end, we first review the $hp$-version finite element method for second order elliptic problems. The optimal convergence of the $hp$-version finite element method with the triangulation containing triangles and parallel elements in two dimension is developed in \cite{BabuvskaSuri1987a} for second order elliptic problems with smooth solutions. We refer to \cite{Schwab1998} for an excellent historical survey. Later on, based on the framework of the Jacobi-weighted Besov spaces (cf. \cite{BabuvskaGuo2001, BabuvskaGuo2002}), the optimal convergence of the $hp$-version finite element method with triangulation containing curvilinear triangles and quadrilaterals in two dimension is established in \cite{GuoSun2007,GuoBabuvska2010} for problems with both smooth and singular solutions.
On the other hand, $hp$-version error estimates are also developed for the DG method of second order problems in recent years (cf.\cite{HoustonSchwabSuli2002, PerugiaSchotzau2002}), which are optimal in the mesh size $h$ and suboptimal in the degree of polynomial $p$.
In \cite{GeorgoulisSuli2005}, under the condition that the exact solution of the reaction-diffusion equation belongs to an augmented Sobolev space, $hp$-optimal error estimates have been deduced for interior penalty DG method with triangulation containing elements being $C^1$-diffeomorphic to parallelograms. And in \cite{StammWihler2010}, by virtue of continuous interpolations of the exact solution, a class of $hp$-version DG methods on parallelograms' mesh for Poisson's equation with homogeneous Dirichlet boundary condition have been proved to converge optimally in the energy norm with respect to both the local element sizes and polynomial degrees.

However, to the best of our knowledge, there are few results about $hp$-version mixed finite element methods for linear elasticity.
Based on the Hellinger-Reissner formulation, a $hp$-mixed finite element space with symmetric stress tensor in two dimensions is constructed in \cite{ArnoldWinther2002}, which is the first stable one using polynomial shape functions.
After establishing the elasticity complex starting from the de Rham complex, Arnold-Falk-Winther element method for the modified Hellinger-Reissner formulation in three space dimensions is devised in \cite{ArnoldFalkWinther2007}, whose stress tensor is nonsymmetric. Furthermore, by designing projection based interpolation operators, Arnold-Falk-Winther element for linear elasticity is extended to meshes with elements of variable order in \cite{QiuDemkowicz2009,QiuDemkowicz2011}.We mention in passing that all the error analyses in these literatures just involve the $h$-version error analysis.

In this paper, we intend to develop the $hp$-version error analysis for the general mixed DG method for the linear elastic problem. To this end, we first derive the $hp$-version error estimates of two $L^2$ projection operators.
Then incorporated with the techniques in \cite{CastilloCockburnPerugiaSchotzau2000} technically, we are able to obtain the $hp$-version error estimates for the previous method in energy norm and $L^2$ norm, respectively. Finally, a variety of numerical examples are provided for validating the theoretical results.

The rest of this paper is organized as follows. Some notations and the DG method in mixed formulation for linear elasticity are
presented in Section 2. The $hp$-version error analysis for the DG method is given in
Section 3. And in Section 4, a series of numerical results are included to show the numerical performance of the DG method proposed.

\section{The DG method for linear elasticity}

Assume that $\Omega\subset \mathbb{R}^d$ $(d=2,3)$ is a bounded
polygon or polyhedron. Let $\boldsymbol{\sigma}=(\sigma_{ij})_{d\times
d}$ be the stress, $\boldsymbol{u}=(u_1,\cdots,u_d)^t$ the
displacement and $\boldsymbol{f}=(f_1,\cdots,f_d)^t$ the applied force. Denote by
$\boldsymbol{\varepsilon}(\boldsymbol{u}):=(\varepsilon_{ij}(\boldsymbol{u}))_{d\times
d}$ the linearized strain tensor with
$\varepsilon_{ij}(\boldsymbol{u}):=(\partial u_i/\partial
x_j+\partial u_j/\partial x_i)/2$, tr the trace operator, and
$\boldsymbol{\div}$ the divergence operator. Consider linear
elasticity in the stress-displacement formulation:
\begin{equation}\label{eq:elas}
\left\{
\begin{array}{ll}
\mathscr{A}\boldsymbol{\sigma}-\boldsymbol{\varepsilon}(\boldsymbol{u})=\boldsymbol{0} & \text{in}\quad \Omega, \\
-\mathbf{div}~\boldsymbol{\sigma}=\boldsymbol{f} & \text{in}\quad \Omega, \\
\boldsymbol{u}=\boldsymbol{0} & \text{on}\quad \partial \Omega,
\end{array}
\right.
\end{equation}
where $\mathscr{A}$ is the compliance tensor of fourth order defined
by
\[
\mathscr{A}\boldsymbol{\sigma}=\frac{1}{2\mu}\left(\boldsymbol{\sigma}-\frac{\lambda}{d\lambda+2\mu}(\textrm{tr}\boldsymbol{\sigma})\boldsymbol{\delta}\right).
\]
Here, $\boldsymbol{\delta}:=(\delta_{ij})_{d\times d}$ is the
Kronecker tensor, and the positive constants $\lambda$ and $\mu$ stand for the
Lam$\acute{e}$ constants.

Then, let us recall the DG method in \cite{ChenHuangHuangXu2010} for solving the problem \eqref{eq:elas}. To this end, we first introduce some notations frequently used later on. For any Banach space $B$, denote by $(B)^s_{d\times d}$ the set of all second order symmetric tensors with entries taking values in $B$. Given a bounded domain
$G\subset\mathbb{R}^{d}$ and a non-negative integer $m$, let
$H^{m}(G)$ be the usual Sobolev space of functions on $G$. The
corresponding norm and semi-norm are denoted respectively by
$\Vert\cdot\Vert_{m,G}$ and $|\cdot|_{m,G}$. If $G$ is $\Omega$, we
abbreviate them by $\Vert\cdot\Vert_{m}$ and $|\cdot|_{m}$,
respectively. Let $H_{0}^{m}(G)$ be the closure of
$C_{0}^{\infty}(G)$ with respect to the norm
$\Vert\cdot\Vert_{m,G}$.

Let $\mathcal {T}_h$ be a regular family of regular triangulations
of $\Omega$ (cf. \cite{Ciarlet1978,BrennerScott2008}). For each $K\in\mathcal
{T}_h$, we denote by $h_K$ the diameter of $K$ and by $\rho_K$ the
diameter of the biggest ball included in $K$. Let $\mathcal{E}_h$ be
the union of all faces of the triangulation $\mathcal {T}_h$ and
$\mathcal{E}^i_h$ be the union of all interior faces of the
triangulation $\mathcal {T}_h$.  The triangulations we consider can have hanging nodes
but have to be regular, that is, there exists a positive constant
$C_1$ such that
\[
\frac{h_K}{\rho_K}\leq C_1 \quad\forall~ K\in \mathcal {T}_h.
\]
Moreover, we restrict the ratio of the sizes of neighbor element
domains. To formally state this property, we need to introduce the
set $\left<K,K'\right>$ defined as follows:
\[
\left<K,K'\right>:=\left\{ \begin{array}{ll} \emptyset, &\textrm{ if
meas}(\partial
K\cap\partial K')=0,\\
\textrm{interior of }\partial K\cap\partial K', &
\textrm{otherwise}.
 \end{array}\right.
\]
Thus we assume that there exists a positive constant $C_2<1$ such
that, for each element $K\in\mathcal {T}_h$,
\[
C_2\leq \frac{h_{K'}}{h_K}\leq \frac{1}{C_2} \quad \forall~
K' \textrm{ satisfying } \left<K,K'\right>\neq \emptyset.
\]
This assumption forbids the situation where the mesh is indefinitely
refined in only one of two adjacent subdomains. Based on the
triangulation $\mathcal{T}_h$, let
\begin{align*}
&\boldsymbol{\Sigma}:=\left\{\boldsymbol{\tau}\in\left(L^2(\Omega)\right)_{d\times
d}^s; \tau_{ij}|_K\in H^1(K)\quad \forall~ K\in\mathcal
{T}_h, i,j=1,\cdots,d\right\}, \\
&\boldsymbol{V}:=\left\{\boldsymbol{v}\in (L^2(\Omega))^d; v_i|_K\in
H^1(K)\quad \forall~ K\in\mathcal {T}_h, i=1,\cdots,d\right\}.
\end{align*}
The corresponding finite element spaces are given by
\begin{align*}
&\boldsymbol{\Sigma}_h:=\left\{\boldsymbol{\tau}\in\left(L^2(\Omega)\right)_{d\times
d}^s; \tau_{ij}|_K\in\mathcal{S}_1(K)\quad \forall~ K\in\mathcal
{T}_h, i,j=1,\cdots,d\right\},\\
&\boldsymbol{V}_h:=\left\{\boldsymbol{v}\in (L^2(\Omega))^d;
v_i|_K\in\mathcal{S}_2(K)\quad \forall~ K\in\mathcal {T}_h,
i=1,\cdots,d\right\},
\end{align*}
where, for each $K\in\mathcal {T}_h$, $\mathcal {S}_1(K)$ and
$\mathcal {S}_2(K)$ are two finite-dimensional spaces of polynomials
in $K$ containing $P_{l_K}(K)$ and $P_{k_K}(K)$, respectively, with integers
$k_K,l_K\geq0$. Here, for a non-negative integer $m$, $P_m(K)$ stands
for the set of all polynomials in $K$ with the total degree no more
than $m$.
We always assume that
\begin{equation}\label{spaces}
\boldsymbol{\varepsilon}(\boldsymbol{V}_h)\subset\boldsymbol{\Sigma}_h, \quad
\boldsymbol{\div}_h\boldsymbol{\Sigma}_h\subset\boldsymbol{V}_h, \quad
\mathscr{A}\boldsymbol{\Sigma}_h\subset\boldsymbol{\Sigma}_h,
\end{equation}
where $\boldsymbol{\div}_h$ is the discrete analogue of the divergence operator
$\boldsymbol{\div}$ with respect to the triangulation $\mathcal {T}_h$, i.e., $(\boldsymbol{\div}_h\boldsymbol{\tau})|_{K}:=\boldsymbol{\div}\boldsymbol{\tau}|_K$
for any $\boldsymbol{\tau}\in\Sigma_h$ and the $K\in \mathcal {T}_h$. It is easy to see from \eqref{spaces} that $|k_K-l_K|\leq 1$ for any $K\in\mathcal {T}_h$.
Then set
\[
p_K=\min\{k_K, l_K\}+1\quad \forall~ K\in\mathcal {T}_h.
\]
Assume that there exists a positive constant $C_3<1$ such
that, for each element $K\in\mathcal {T}_h$,
\[
C_3\leq \frac{p_{K'}}{p_K}\leq \frac{1}{C_3} \quad \forall~
K' \textrm{ satisfying } \left<K,K'\right>\neq \emptyset.
\]
For a function $v\in
L^{2}(\Omega)$ with $v|_{K}\in H^{m}(K)$ for all
$K\in\mathcal{T}_{h}$, let $\Vert v\Vert_{m,h}$ and $|v|_{m,h}$ be
the usual broken $H^m$-type norm and semi-norm of $v$:
\[
\Vert v\Vert_{m,h}:=\bigg(\sum_{K\in\mathcal{T}_{h}}\Vert
v\Vert_{m,K}^{2}\bigg)^{1/2},\;\;|v|_{m,h}:=\bigg(\sum_{K\in
\mathcal{T}_{h}}|v|_{m,K}^{2}\bigg)^{1/2}.
\]
If $v$ is a vector-value or tensor-value function, the corresponding
$\|\cdot\|_{m,h}$ and $|\cdot|_{m,h}$ are defined in the similar
manners. For a vector or tensor $\boldsymbol{v}$, its length
$|\boldsymbol{v}|$ is $(\boldsymbol{v}\cdot\boldsymbol{v})^{1/2}$ or
$(\boldsymbol{v}:\boldsymbol{v})^{1/2}$. Here the symbol $:$ denotes
the double dot product operation of tensors. Throughout this paper,
we use the notation \textquotedblleft\ $\lesssim\cdots
$\textquotedblright\ to mean that \textquotedblleft$\leq C\cdots
$\textquotedblright, where $C$ is a generic positive constant
independent of local element sizes and polynomial degrees, which may take different
values at different appearances. And $a\eqsim b$ means $a\lesssim b$ and $b\lesssim a$.

Let $K^+$ and $K^-$ be two adjacent elements of $\mathcal{T}_h$. Let
$\boldsymbol{x}$ be an arbitrary point of the set
$e'=\left<K^+,K^-\right>$, and let $\boldsymbol{n}^+$ and
$\boldsymbol{n}^-$ be the corresponding outward unit normals at that
point. For a vector-valued function $\boldsymbol{v}$ and tensor-valued function $\boldsymbol{\tau}$ smooth inside
each element $K^\pm$, let us denote by $\boldsymbol{v}^\pm$ and $\boldsymbol{\tau}^\pm$ the
trace of $\boldsymbol{v}$ and $\boldsymbol{\tau}$ on $e'$ from the interior of $K^\pm$, respectively. Then
we define averages and jumps at $\boldsymbol{x}\in e'$ as follows:
\begin{align*}
\{\boldsymbol{v}\}=\frac{1}{2}(\boldsymbol{v}^++\boldsymbol{v}^-),
&\quad
[\boldsymbol{v}]=\boldsymbol{v}^+\cdot\boldsymbol{n}^++\boldsymbol{v}^-\cdot\boldsymbol{n}^-,
\\
\{\boldsymbol{\tau}\}=\frac{1}{2}(\boldsymbol{\tau}^++\boldsymbol{\tau}^-),
&\quad
[\boldsymbol{\tau}]=\boldsymbol{\tau}^+\boldsymbol{n}^++\boldsymbol{\tau}^-\boldsymbol{n}^-.
\end{align*}
If $\boldsymbol{x}$ is on an face $e$ lying on the boundary
$\partial\Omega$, the above terms are defined by
\begin{align*}
\{\boldsymbol{v}\}=\boldsymbol{v}, &\quad
[\boldsymbol{v}]=\boldsymbol{v}\cdot\boldsymbol{n},
\\
\{\boldsymbol{\tau}\}=\boldsymbol{\tau}, &\quad
[\boldsymbol{\tau}]=\boldsymbol{\tau}\boldsymbol{n},
\end{align*}
where $\boldsymbol{n}$ is the unit outward normal vector on
$\partial\Omega$. In addition, we define a matrix valued jump
$\llbracket\cdot\rrbracket$ of a vector $\boldsymbol{v}$ as follows:
\begin{align*}
&\llbracket\boldsymbol{v}\rrbracket=\frac{1}{2}(\boldsymbol{v}^+\otimes\boldsymbol{n}^++\boldsymbol{n}^+\otimes\boldsymbol{v}^++\boldsymbol{v}^-\otimes\boldsymbol{n}^-+\boldsymbol{n}^-\otimes\boldsymbol{v}^-),
\quad\textrm{if} \; \boldsymbol{x}\in e \in \mathcal{E}_h^i,\\
&\llbracket\boldsymbol{v}\rrbracket=\frac{1}{2}(\boldsymbol{v}\otimes\boldsymbol{n}+\boldsymbol{n}\otimes\boldsymbol{v}),
\quad \textrm{if} \; \boldsymbol{x}\in e \in
\mathcal{E}_h\cap\partial\Omega,
\end{align*}
where $\boldsymbol{v}\otimes\boldsymbol{n}$ denote the matrix whose
$(i,j)$-th entry is $v_in_j$ for two vectors $\boldsymbol{v}$ and
$\boldsymbol{n}$.

With the help of the above notations, the mixed DG method devised in \cite{ChenHuangHuangXu2010} for linear elasticity problem~\eqref{eq:elas} can be described as follows.

Find
$(\boldsymbol{\sigma}_h,\boldsymbol{u}_h)\in\boldsymbol{\Sigma}_h\times\boldsymbol{V}_h$
such that
\begin{align}
 a(\boldsymbol{\sigma}_h,\boldsymbol{\tau})+b(\boldsymbol{u}_h,\boldsymbol{\tau})&=0,  \label{eq:mix1}\\
-b(\boldsymbol{v},\boldsymbol{\sigma}_h)+c(\boldsymbol{u}_h,\boldsymbol{v})&=F(\boldsymbol{v}),\label{eq:mix2}
\end{align}
for all
$(\boldsymbol{\tau},\boldsymbol{v})\in\boldsymbol{\Sigma}_h\times\boldsymbol{V}_h$,
where
\begin{align*}
&a(\boldsymbol{\sigma},\boldsymbol{\tau}):=\int_\Omega\mathscr{A}\boldsymbol{\sigma}:\boldsymbol{\tau}dx+\int_{\mathcal{E}^i_{h}}C_{22}[\boldsymbol{\sigma}]\cdot[\boldsymbol{\tau}]ds,
\\
&
b(\boldsymbol{v},\boldsymbol{\tau}):=-\sum_{K\in\mathcal{T}_h}\int_K\boldsymbol{\varepsilon}(\boldsymbol{v}):\boldsymbol{\tau}dx+\int_{\mathcal{E}_{h}}\llbracket\boldsymbol{v}\rrbracket:\{\boldsymbol{\tau}\}ds,
\\
&c(\boldsymbol{u},\boldsymbol{v}):=\int_{\mathcal{E}_h}C_{11}\llbracket\boldsymbol{u}\rrbracket:\llbracket\boldsymbol{v}\rrbracket
ds,
\\
& F(\boldsymbol{v}):=\int_\Omega\boldsymbol{f}\cdot\boldsymbol{v}dx,
\end{align*}
with $C_{11}>0, C_{22}\geq 0$.

\begin{remark}\rm
 If $C_{22}= 0$, the above method is reduced to the LDG method in \cite{ChenHuangHuangXu2010}, whose $h$-version error analysis has also been well studied  there. Here, we will focus on $hp$-version error estimates for the method for $C_{22}> 0$. It deserves to point out that our arguments developed in this paper can not applies to the case where $C_{22}=0$, since the mesh-dependent norm $|\cdot|_B$ used in the next section does not make sense in this case.
 \end{remark}

\begin{remark}\rm
For simplicity, we confine ourselves to error analysis for the mixed DG method \eqref{eq:mix1}-\eqref{eq:mix2} related to homogenous Dirichlet boundary conditions. As a matter of fact, the mathematical reasoning developed in what follows can be naturally extended to deal with a mixed DG method for the elastic problem with general mixed boundary conditions:
\begin{equation*}
\left\{
\begin{array}{ll}
\mathscr{A}\boldsymbol{\sigma}-\boldsymbol{\varepsilon}(\boldsymbol{u})=\boldsymbol{0} & \text{in}\quad \Omega, \\
-\mathbf{div}~\boldsymbol{\sigma}=\boldsymbol{f} & \text{in}\quad \Omega, \\
\boldsymbol{u}=\boldsymbol{g}_D & \text{on}\quad \Gamma_D, \\
\boldsymbol{\sigma}\boldsymbol{n}=\boldsymbol{g}_N\boldsymbol{n} & \text{on}\quad \Gamma_N,
\end{array}
\right.
\end{equation*}
where $\Gamma_D$ and $\Gamma_N$ are two disjoint subsets of $\partial\Omega$ such that $\mbox{meas}(\Gamma_D)\ne 0$ and $\overline{\Gamma_D\cup\Gamma_N}=\partial\Omega$. In this case, following the ideas in \cite{ChenHuangHuangXu2010},
the corresponding mixed DG method is to find
$(\boldsymbol{\sigma}_h,\boldsymbol{u}_h)\in\boldsymbol{\Sigma}_h\times\boldsymbol{V}_h$
such that
\begin{align}
 \tilde a(\boldsymbol{\sigma}_h,\boldsymbol{\tau})+\tilde b(\boldsymbol{u}_h,\boldsymbol{\tau})&=\tilde F_1(\boldsymbol{\tau}),  \label{eq:mixDN1}\\
-\tilde b(\boldsymbol{v},\boldsymbol{\sigma}_h)+\tilde c(\boldsymbol{u}_h,\boldsymbol{v})&=\tilde F_2(\boldsymbol{v}),\label{eq:mixDN2}
\end{align}
for all
$(\boldsymbol{\tau},\boldsymbol{v})\in\boldsymbol{\Sigma}_h\times\boldsymbol{V}_h$,
where
\begin{align*}
&\tilde a(\boldsymbol{\sigma},\boldsymbol{\tau}):=\int_\Omega\mathscr{A}\boldsymbol{\sigma}:\boldsymbol{\tau}dx+\int_{\mathcal{E}^i_{h}}C_{22}[\boldsymbol{\sigma}]\cdot[\boldsymbol{\tau}]ds+\int_{\Gamma_N}C_{22}(\boldsymbol{\sigma}\boldsymbol{n})\cdot(\boldsymbol{\tau}\boldsymbol{n})ds,
\\
&
\tilde b(\boldsymbol{v},\boldsymbol{\tau}):=-\sum_{K\in\mathcal{T}_h}\int_K\boldsymbol{\varepsilon}(\boldsymbol{v}):\boldsymbol{\tau}dx+\int_{\mathcal{E}_{h}^i}\llbracket\boldsymbol{v}\rrbracket:\{\boldsymbol{\tau}\}ds +\int_{\Gamma_D}\boldsymbol{v}\cdot(\boldsymbol{\tau}\boldsymbol{n})ds,
\\
&\tilde c(\boldsymbol{u},\boldsymbol{v}):=\int_{\mathcal{E}_h^i}C_{11}\llbracket\boldsymbol{u}\rrbracket:\llbracket\boldsymbol{v}\rrbracket
ds + \int_{\Gamma_D}C_{11}\llbracket\boldsymbol{u}\rrbracket:\llbracket\boldsymbol{v}\rrbracket
ds,
\\
& \tilde F_1(\boldsymbol{\tau}):=\int_{\Gamma_D}\boldsymbol{g}_D\cdot(\boldsymbol{\tau}\boldsymbol{n})ds+\int_{\Gamma_N}C_{22}(\boldsymbol{g}_N\boldsymbol{n})\cdot(\boldsymbol{\tau}\boldsymbol{n})ds,
\\
& \tilde F_2(\boldsymbol{v}):=\int_\Omega\boldsymbol{f}\cdot\boldsymbol{v}dx + \int_{\Gamma_D}C_{11}\llbracket\boldsymbol{g}_D\rrbracket:\llbracket\boldsymbol{v}\rrbracket ds +\int_{\Gamma_N}(\boldsymbol{g}_N\boldsymbol{n})\cdot\boldsymbol{v}ds.
\end{align*}
\end{remark}

\section{The $hp$-version error analysis for the DG method}

In this section, we are going to establish $hp$-version error estimates for the DG method
\eqref{eq:mix1}-\eqref{eq:mix2}. Our derivation is mainly based on the techniques developed in \cite{CastilloCockburnPerugiaSchotzau2000}. To this end, we first rewrite
\eqref{eq:mix1}-\eqref{eq:mix2} in a compact form, described as follows.

Find
$(\boldsymbol{\sigma}_h,\boldsymbol{u}_h)\in\boldsymbol{\Sigma}_h\times\boldsymbol{V}_h$
such that
\begin{equation}\label{eq:variation}
A(\boldsymbol{\sigma}_h,\boldsymbol{u}_h;\boldsymbol{\tau},\boldsymbol{v})=F(\boldsymbol{v})
\end{equation}
for all
$(\boldsymbol{\tau},\boldsymbol{v})\in\boldsymbol{\Sigma}_h\times\boldsymbol{V}_h$,
where
\begin{equation}
\label{definition of A}
A(\boldsymbol{\sigma},\boldsymbol{u};\boldsymbol{\tau},\boldsymbol{v}):=a(\boldsymbol{\sigma},\boldsymbol{\tau})+b(\boldsymbol{u},\boldsymbol{\tau})-b(\boldsymbol{v},\boldsymbol{\sigma})+c(\boldsymbol{u},\boldsymbol{v}).
\end{equation}
In the following, we always assume that
$(\boldsymbol{\sigma},\boldsymbol{u})\in
\left(H^{1}(\Omega)\right)_{d\times d}^s\times (H^{2}(\Omega)\cap H_0^{1}(\Omega))^d$ is the solution of the original
problem \eqref{eq:elas}. Let $\boldsymbol{P}_h$ be $L^2$ projection
operator from $\boldsymbol{\Sigma}$ onto the finite element space
$\boldsymbol{\Sigma}_h$ and $\boldsymbol{Q}_h$ be $L^2$ projection
operator from $\boldsymbol{V}$ onto the finite element space
$\boldsymbol{V}_h$. For simplicity, we still write
$\boldsymbol{P}_h$ and $\boldsymbol{Q}_h$ for $\boldsymbol{P}_h|_K$
and $\boldsymbol{Q}_h|_K$.

From Theorem~1.1 in \cite{Chernov2012}, Corollary~1.2 in \cite{MelenkWurzer2014}, Lemma~4.4 in \cite{BabuvskaSuri1987a} and the scaling argument,
we can easily obtain the following $hp$-version error estimates for $L^2$ projection
operators $\boldsymbol{P}_h$ and $\boldsymbol{Q}_h$.

\begin{lem}\label{lem:interpolation}
Let $\boldsymbol{\tau}\in \left(H^{s+1}(K)\right)_{d\times d}^s$,
$\boldsymbol{v}\in H^{s+2}(K)^{d}$, $s\geq0$. Then
\begin{align*}
\left\|\boldsymbol{\tau}-\boldsymbol{P}_h\boldsymbol{\tau}\right\|_{0,K}+\frac{h_K}{l_K+1}\left|\boldsymbol{\tau}-\boldsymbol{P}_h\boldsymbol{\tau}\right|_{1,K}&\lesssim
\left(\frac{h_K}{l_K+1}\right)^{r+1} \|\boldsymbol{\tau}\|_{r+1,K}, \\
\left\|\boldsymbol{\tau}-\boldsymbol{P}_h\boldsymbol{\tau}\right\|_{0,\partial
K}&\lesssim
\left(\frac{h_K}{l_K+1}\right)^{r+1/2} \|\boldsymbol{\tau}\|_{r+1,K},
\end{align*}
for $0\leq r\leq \min\{s,l_K\}$, and
\begin{align*}
\left\|\boldsymbol{v}-\boldsymbol{Q}_h\boldsymbol{v}\right\|_{0,K}+\frac{h_K}{k_K+1}\left|\boldsymbol{v}-\boldsymbol{Q}_h\boldsymbol{v}\right|_{1,K}&\lesssim
\left(\frac{h_K}{k_K+1}\right)^{r+1}\|\boldsymbol{v}\|_{r+1,K}, \\
\left\|\boldsymbol{v}-\boldsymbol{Q}_h\boldsymbol{v}\right\|_{0,\partial
K}&\lesssim
\left(\frac{h_K}{k_K+1}\right)^{r+1/2}\|\boldsymbol{v}\|_{r+1,K},
\end{align*}
for $0\leq r\leq \min\{s+1,k_K\}$.
\end{lem}\vspace{1 pc}

The next result shows the Galerkin orthogonality holds true for the numerical solution of the method \eqref{eq:variation} (or equivalently, the method \eqref{eq:mix1}-\eqref{eq:mix2}).
\begin{lem}
Let $(\boldsymbol{\sigma},\boldsymbol{u})$ be the solution of problem \eqref{eq:elas}, and let $(\boldsymbol{\sigma}_h,\boldsymbol{u}_h)$ be the solution of the DG method \eqref{eq:variation}. Then for any $(\boldsymbol{\tau},\boldsymbol{v})\in\boldsymbol{\Sigma}_h\times\boldsymbol{V}_h$,
there holds
\begin{equation}\label{eq:go}
A(\boldsymbol{\sigma}-\boldsymbol{\sigma}_h,\boldsymbol{u}-\boldsymbol{u}_h;\boldsymbol{\tau},\boldsymbol{v})=0.
\end{equation}
\end{lem}
\noindent\emph{Proof}. By the definition of $A$, we have
\begin{align*}
A(\boldsymbol{\sigma},\boldsymbol{u};\boldsymbol{\tau},\boldsymbol{v})=&\int_{\Omega}\mathscr{A}\boldsymbol{\sigma}:\boldsymbol{\tau}dx+\int_{\mathcal{E}_h^i}C_{22}[\boldsymbol{\sigma}]\cdot[\boldsymbol{\tau}]ds-\sum_{K\in\mathcal{T}_h}\int_K\boldsymbol{\varepsilon}(\boldsymbol{u}):\boldsymbol{\tau}dx \\
&+\int_{\mathcal{E}_{h}}\llbracket\boldsymbol{u}\rrbracket:\{\boldsymbol{\tau}\}ds +\sum_{K\in\mathcal{T}_h}\int_K\boldsymbol{\varepsilon}(\boldsymbol{v}):\boldsymbol{\sigma}dx-\int_{\mathcal{E}_{h}}\llbracket\boldsymbol{v}\rrbracket:\{\boldsymbol{\sigma}\}ds
\\
&+\int_{\mathcal{E}_h}C_{11}\llbracket\boldsymbol{u}\rrbracket:\llbracket\boldsymbol{v}\rrbracket
ds.
\end{align*}
Since $(\boldsymbol{\sigma},\boldsymbol{u})\in
\left(H^{1}(\Omega)\right)_{d\times d}^s\times (H^{2}(\Omega)\cap H_0^{1}(\Omega))^d$ is
the solution of problem \eqref{eq:elas},
the quantities $\llbracket\boldsymbol{u}\rrbracket$ and $[\boldsymbol{\sigma}]$
both vanish. Hence, we can rewrite
$A(\boldsymbol{\sigma},\boldsymbol{u};\boldsymbol{\tau},\boldsymbol{v})$
as
\begin{align*}
A(\boldsymbol{\sigma},\boldsymbol{u};\boldsymbol{\tau},\boldsymbol{v})=&\int_{\Omega}\mathscr{A}\boldsymbol{\sigma}:\boldsymbol{\tau}dx-\sum_{K\in\mathcal{T}_h}\int_K\boldsymbol{\varepsilon}(\boldsymbol{u}):\boldsymbol{\tau}dx
+\sum_{K\in\mathcal{T}_h}\int_K\boldsymbol{\varepsilon}(\boldsymbol{v}):\boldsymbol{\sigma}dx-\int_{\mathcal{E}_{h}}\llbracket\boldsymbol{v}\rrbracket:\{\boldsymbol{\sigma}\}ds.
\end{align*}
By \eqref{eq:elas} and integration by parts, we then have
\begin{equation}\label{eq:variation1}
A(\boldsymbol{\sigma},\boldsymbol{u};\boldsymbol{\tau},\boldsymbol{v})=F(\boldsymbol{v}),
\end{equation}
from which and \eqref{eq:variation}, the desired identity \eqref{eq:go} follows readily. \quad\quad $\Box$\vspace{1 pc}

To derive our error analysis, we still require to establish a number of inequalities revealing the approximation properties of the projection operators $\boldsymbol{P}_h$ and $\boldsymbol{Q}_h$. Before doing this, we first introduce two
seminorms for later requirement. For
$(\boldsymbol{\tau},\boldsymbol{v})\in\boldsymbol{\Sigma}\times\boldsymbol{V}$,
define
\begin{align*}
|(\boldsymbol{\tau},\boldsymbol{v})|^2_A &=A(\boldsymbol{\tau},\boldsymbol{v};\boldsymbol{\tau},\boldsymbol{v}) =\int_{\Omega}\mathscr{A}\boldsymbol{\tau}:\boldsymbol{\tau}dx+\int_{\mathcal
{E}^i_{h}}\left(C_{22}[\boldsymbol{\tau}]^2+C_{11}\llbracket\boldsymbol{v}\rrbracket^2\right)ds
+\int_{\partial\Omega}C_{11}\llbracket\boldsymbol{v}\rrbracket^2ds,
\\
|(\boldsymbol{\tau},\boldsymbol{v})|^2_B &=\int_{\mathcal
{E}^i_{h}}\left(C_{22}[\boldsymbol{\tau}]^2+\frac{1}{C_{11}}\{\boldsymbol{\tau}\}^2+\frac{1}{C_{22}}\{\boldsymbol{v}\}^2+C_{11}\llbracket\boldsymbol{v}\rrbracket^2\right)ds
+\int_{\partial\Omega}\left(\frac{1}{C_{11}}|\boldsymbol{\tau}|^2+C_{11}\llbracket\boldsymbol{v}\rrbracket^2\right)ds.
\end{align*}

And we also want to introduce two functionals $K_A$ and $K_B$, by which all the
error estimates we are interested in can be obtained. For
$(\boldsymbol{\sigma},\boldsymbol{u})\in (H^{s+1}(\Omega))_{d\times
d}^s\times H^{s+2}(\Omega)^d$ and
$(\boldsymbol{\tau},\boldsymbol{v})\in (H^{t+1}(\Omega))_{d\times
d}^s\times H^{t+2}(\Omega)^d$ with $s,t\geq0$, define
\[
K_A(\boldsymbol{\sigma},\boldsymbol{u};\boldsymbol{\tau},\boldsymbol{v})=\left\{
\begin{array}{ll}
\sum\limits^5_{i=1}S_i(\boldsymbol{\sigma},\boldsymbol{u};\boldsymbol{\tau},\boldsymbol{v}),
& \textrm{if }  (\boldsymbol{\sigma},\boldsymbol{u})\neq(\boldsymbol{\tau},\boldsymbol{v}),\\
\sum\limits_{i=1,2,5}S_i(\boldsymbol{\sigma},\boldsymbol{u};\boldsymbol{\sigma},\boldsymbol{u}),& \textrm{if } (\boldsymbol{\sigma},\boldsymbol{u})=(\boldsymbol{\tau},\boldsymbol{v}),\\
\end{array}
\right.
\]
where
\begin{align*}
S_1=&\left(\sum_{K\in\mathcal{T}_h}\left(\frac{h_K}{p_K}\right)^{2s_{1K}+2}\|\boldsymbol{\sigma}\|^2_{s_{1K}+1,K}\right)^{\frac{1}{2}}\left(\sum_{K\in\mathcal{T}_h}\left(\frac{h_K}{p_K}\right)^{2t_{1K}+2}\|\boldsymbol{\tau}\|^2_{t_{1K}+1,K}\right)^{\frac{1}{2}},
\end{align*}
\begin{align*}
S_2=&\left(\sum_{K\in\mathcal{T}_h}C^{\partial
K}_{22}\left(\frac{h_K}{p_K}\right)^{2s_{1K}+1}\|\boldsymbol{\sigma}\|^2_{s_{1K}+1,K}\right)^{\frac{1}{2}} \left(\sum_{K\in\mathcal{T}_h}C^{\partial
K}_{22}\left(\frac{h_K}{p_K}\right)^{2t_{1K}+1}\|\boldsymbol{\tau}\|^2_{t_{1K}+1,K}\right)^{\frac{1}{2}},
\end{align*}
\begin{align*}
S_3=&\left(\sum_{K\in\mathcal{T}_h}C^{\partial
K}_{11}\left(\frac{h_K}{p_K}\right)^{2s_{2K}+1}\|\boldsymbol{u}\|^2_{s_{2K}+1,K}\right)^{\frac{1}{2}} \left(\sum_{K\in\mathcal{T}_h}\frac{1}{\widetilde{C}^{\partial
K}_{11}}\left(\frac{h_K}{p_K}\right)^{2t_{1K}+1}\|\boldsymbol{\tau}\|^2_{t_{1K}+1,K}\right)^{\frac{1}{2}},
\end{align*}
\begin{align*}
S_4=&\left(\sum_{K\in\mathcal{T}_h}\frac{1}{\widetilde{C}^{\partial
K}_{11}}\left(\frac{h_K}{p_K}\right)^{2s_{1K}+1}\|\boldsymbol{\sigma}\|^2_{s_{1K}+1,K}\right)^{\frac{1}{2}} \left(\sum_{K\in\mathcal{T}_h}C^{\partial
K}_{11}\left(\frac{h_K}{p_K}\right)^{2t_{2K}+1}\|\boldsymbol{v}\|^2_{t_{2K}+1,K}\right)^{\frac{1}{2}},
\end{align*}
\begin{align*}
S_5=&\left(\sum_{K\in\mathcal{T}_h}C^{\partial
K}_{11}\left(\frac{h_K}{p_K}\right)^{2s_{2K}+1}\|\boldsymbol{u}\|^2_{s_{2K}+1,K}\right)^{\frac{1}{2}} \left(\sum_{K\in\mathcal{T}_h}C^{\partial
K}_{11}\left(\frac{h_K}{p_K}\right)^{2t_{2K}+1}\|\boldsymbol{v}\|^2_{t_{2K}+1,K}\right)^{\frac{1}{2}}
\end{align*}
with $0\leq s_{1K}\leq \min\{s,l_K\}$, $0\leq s_{2K}\leq \min\{s+1,k_K\}$, $0\leq t_{1K}\leq \min\{t,l_K\}$, $0\leq t_{2K}\leq \min\{t+1,k_K\}$,
 $\widetilde{C}^{\partial
K}_{ii}:=\inf\{C_{ii}(\boldsymbol{x}); \boldsymbol{x}\in\partial K\}$, $C^{\partial
K}_{ii}:=\sup\{C_{ii}(\boldsymbol{x}); \boldsymbol{x}\in\partial K\}$ for
$i=1,2$.
The quantity $K_B$ is defined as
\begin{align*}
K^2_B(\boldsymbol{\sigma},\boldsymbol{u})=&\sum_{K\in\mathcal{T}_h}\left(\left(\frac{h_K}{p_K}\right)^{2s_{1K}+1}\left(\frac{1}{\widetilde{C}^{\partial
K}_{11}}+C^{\partial
K}_{22}\right)\|\boldsymbol{\sigma}\|^2_{s_{1K}+1,K}\right) \\
&+\sum_{K\in\mathcal{T}_h}\left(\left(\frac{h_K}{p_K}\right)^{2s_{2K}+1}\left(C^{\partial
K}_{11}+\frac{1}{\widetilde{C}_{22}^{\partial
K}}\right)\|\boldsymbol{u}\|^2_{s_{2K}+1,K}\right).
\end{align*}


\begin{lem}\label{lem:1}
For any
$(\boldsymbol{\sigma},\boldsymbol{u}),(\boldsymbol{\tau},\boldsymbol{v})\in
\boldsymbol{\Sigma}\times \boldsymbol{V}$, assume that for each
$K\in\mathcal {T}_h$,
$(\boldsymbol{\sigma},\boldsymbol{u})|_K\in\left(H^{s+1}(K)\right)_{d\times
d}^s\times H^{s+2}(K)^d$ and
$(\boldsymbol{\tau},\boldsymbol{v})|_K\in\left(H^{t+1}(K)\right)_{d\times
d}^s\times H^{t+2}(K)^d$. Then
\[
A(\boldsymbol{\sigma}-\boldsymbol{P}_h\boldsymbol{\sigma},\boldsymbol{u}-\boldsymbol{Q}_h\boldsymbol{u};\boldsymbol{\tau}-\boldsymbol{P}_h\boldsymbol{\tau},\boldsymbol{v}-\boldsymbol{Q}_h\boldsymbol{v})\lesssim
K_A(\boldsymbol{\sigma},\boldsymbol{u};\boldsymbol{\tau},\boldsymbol{v}).
\]
\end{lem}

\noindent\emph{Proof}. For convenience, set
$\boldsymbol{\xi_{\sigma}}:=\boldsymbol{\sigma}-\boldsymbol{P}_h\boldsymbol{\sigma},~
\boldsymbol{\xi_{u}}:=\boldsymbol{u}-\boldsymbol{Q}_h\boldsymbol{u},~
\boldsymbol{\xi_{\tau}}:=\boldsymbol{\tau}-\boldsymbol{P}_h\boldsymbol{\tau},~
\boldsymbol{\xi_{v}}:=\boldsymbol{v}-\boldsymbol{Q}_h\boldsymbol{v}$.
We start by writing
\[
A(\boldsymbol{\xi_{\sigma}},\boldsymbol{\xi_{u}};\boldsymbol{\xi_{\tau}},\boldsymbol{\xi_{v}}):=a(\boldsymbol{\xi_{\sigma}},\boldsymbol{\xi_{\tau}})+b(\boldsymbol{\xi_{u}},\boldsymbol{\xi_{\tau}})-b(\boldsymbol{\xi_{v}},\boldsymbol{\xi_{\sigma}})+c(\boldsymbol{\xi_{u}},\boldsymbol{\xi_{v}}),
\]
and then proceed by estimating each term on the right-hand side
separately. According to the Cauchy-Schwarz inequality and Lemma
\ref{lem:interpolation}, we have
\begin{align*}
\left|a(\boldsymbol{\xi_{\sigma}},\boldsymbol{\xi_{\tau}})\right|&=\left|\sum_{K\in\mathcal{T}_h}\int_K\mathcal{A}\boldsymbol{\xi_{\sigma}}:\boldsymbol{\xi_{\tau}}dx+\sum_{e\in\mathcal{E}^i_{h}}\int_{e}C_{22}[\boldsymbol{\xi_{\sigma}}]\cdot[\boldsymbol{\xi_{\tau}}]ds\right|
\\
& \lesssim
\sum_{K\in\mathcal{T}_h}\|\boldsymbol{\xi_{\sigma}}\|_{0,K}\|\boldsymbol{\xi_{\tau}}\|_{0,K}+\sum_{e\in\mathcal{E}^i_{h}}\|\sqrt{C_{22}}[\boldsymbol{\xi_{\sigma}}]\|_{0,e}\|\sqrt{C_{22}}[\boldsymbol{\xi_{\tau}}]\|_{0,e}
\\
& \lesssim
\left(\sum_{K\in\mathcal{T}_h}\|\boldsymbol{\xi_{\sigma}}\|^2_{0,K}\right)^{\frac{1}{2}}\left(\sum_{K\in\mathcal{T}_h}\|\boldsymbol{\xi_{\tau}}\|^2_{0,K}\right)^{\frac{1}{2}} +\left(\sum_{e\in\mathcal{E}^i_{h}}\|\sqrt{C_{22}}[\boldsymbol{\xi_{\sigma}}]\|^2_{0,e}\right)^{\frac{1}{2}}\left(\sum_{e\in\mathcal{E}^i_{h}}\|\sqrt{C_{22}}[\boldsymbol{\xi_{\tau}}]\|^2_{0,e}\right)^{\frac{1}{2}}
\\
&\lesssim
S_1(\boldsymbol{\sigma},\boldsymbol{u};\boldsymbol{\tau},\boldsymbol{v})+S_2(\boldsymbol{\sigma},\boldsymbol{u};\boldsymbol{\tau},\boldsymbol{v}).
\end{align*}
Again, by the Cauchy-Schwarz  inequality and
Lemma \ref{lem:interpolation}, it follows that
\begin{align*}
\left|b(\boldsymbol{\xi_{u}},\boldsymbol{\xi_{\tau}})\right|&=\left|-\sum_{K\in\mathcal{T}_h}\int_K\boldsymbol{\varepsilon}(\boldsymbol{\xi_{u}}):\boldsymbol{\xi_{\tau}}dx+\int_{\mathcal{E}_{h}}\llbracket\boldsymbol{\xi_{u}}\rrbracket:\{\boldsymbol{\xi_{\tau}}\}ds\right|
\\
&\lesssim
\left(\sum_{K\in\mathcal{T}_h}C^{\partial
K}_{11}\left(\frac{h_K}{p_K}\|\nabla\boldsymbol{\xi_{u}}\|^2_{0,K}+\|\boldsymbol{\xi_{u}}\|^2_{0,\partial
K}\right)\right)^{\frac{1}{2}}  \cdot\left(\sum_{K\in\mathcal{T}_h}\frac{1}{\widetilde{C}^{\partial
K}_{11}}\left(\frac{p_K}{h_K}\|\boldsymbol{\xi_{\tau}}\|^2_{0,K}+\|\boldsymbol{\xi_{\tau}}\|^2_{0,\partial
K}\right)\right)^{\frac{1}{2}} \\
&\lesssim
S_3(\boldsymbol{\sigma},\boldsymbol{u};\boldsymbol{\tau},\boldsymbol{v}).
\end{align*}
Using the similar arguments, we can also derive
\begin{align*}
\left|b(\boldsymbol{\xi_{v}},\boldsymbol{\xi_{\sigma}})\right|&\lesssim
S_4(\boldsymbol{\sigma},\boldsymbol{u};\boldsymbol{\tau},\boldsymbol{v}),\\
\left|c(\boldsymbol{\xi_{u}},\boldsymbol{\xi_{v}})\right|&\lesssim
\left(\sum_{K\in\mathcal{T}_h}C^{\partial
K}_{11}\|\boldsymbol{\xi_{u}}\|^2_{0,\partial
K}\right)^{\frac{1}{2}}\left(\sum_{K\in\mathcal{T}_h}C^{\partial
K}_{11}\|\boldsymbol{\xi_{v}}\|^2_{0,\partial K}\right)^{\frac{1}{2}} \lesssim
S_5(\boldsymbol{\sigma},\boldsymbol{u};\boldsymbol{\tau},\boldsymbol{v}).
\end{align*}
This proves the required estimate for
$(\boldsymbol{\sigma},\boldsymbol{u})\neq(\boldsymbol{\tau},\boldsymbol{v})$.
If $(\boldsymbol{\sigma},\boldsymbol{u})=(\boldsymbol{\tau},\boldsymbol{v})$, the required estimate follows immediately from the identity
\[
A(\boldsymbol{\xi_{\sigma}},\boldsymbol{\xi_{u}};\boldsymbol{\xi_{\sigma}},\boldsymbol{\xi_{u}})=a(\boldsymbol{\xi_{\sigma}},\boldsymbol{\xi_{\sigma}})+c(\boldsymbol{\xi_{u}},\boldsymbol{\xi_{u}}).
\quad\quad\quad\Box
\]

\begin{lem}\label{lem:2}
For any $(\boldsymbol{\sigma},\boldsymbol{u})\in
\boldsymbol{\Sigma}\times
\boldsymbol{V},(\boldsymbol{\tau},\boldsymbol{v})\in
\boldsymbol{\Sigma}_h\times \boldsymbol{V}_h$, there holds
\[
A(\boldsymbol{\tau},\boldsymbol{v};\boldsymbol{\sigma}-\boldsymbol{P}_h\boldsymbol{\sigma},\boldsymbol{u}-\boldsymbol{Q}_h\boldsymbol{u})\lesssim
\left|(\boldsymbol{\tau},\boldsymbol{v})\right|_A\left|(\boldsymbol{\sigma}-\boldsymbol{P}_h\boldsymbol{\sigma},\boldsymbol{u}-\boldsymbol{Q}_h\boldsymbol{u})\right|_B.
\]
\end{lem}

\noindent\emph{Proof}. By setting
$\boldsymbol{\xi_{\sigma}}:=\boldsymbol{\sigma}-\boldsymbol{P}_h\boldsymbol{\sigma},
\boldsymbol{\xi_{u}}:=\boldsymbol{u}-\boldsymbol{Q}_h\boldsymbol{u}$,
we have
\begin{align*}
\left|A(\boldsymbol{\tau},\boldsymbol{v};\boldsymbol{\xi_{\sigma}},\boldsymbol{\xi_{u}})\right|
&
\leq\left|a(\boldsymbol{\tau},\boldsymbol{\xi_{\sigma}})\right|+\left|b(\boldsymbol{v},\boldsymbol{\xi_{\sigma}})\right|+\left|b(\boldsymbol{\xi_{u}},\boldsymbol{\tau})\right|+\left|c(\boldsymbol{v},\boldsymbol{\xi_{u}})\right|
\\
&=:T_1+T_2+T_3+T_4.
\end{align*}
By the inclusion property \eqref{spaces}, we have
$\int_K\mathscr{A}\boldsymbol{\tau}:\boldsymbol{\xi_{\tau}}dx=0$.
Hence, we have by the Cauchy-Schwarz inequality that
\[
T_1\leq
\left(\int_{\mathcal{E}^i_{h}}C_{22}[\boldsymbol{\tau}]^2ds\right)^{\frac{1}{2}}\left(\int_{\mathcal{E}^i_{h}}C_{22}[\boldsymbol{\xi_{\sigma}}]^2ds\right)^{\frac{1}{2}}\leq\left|(\boldsymbol{\tau},\boldsymbol{v})\right|_A\left|(\boldsymbol{\xi_{\sigma}},\boldsymbol{\xi_{u}})\right|_B.
\]
Furthermore, since
$\int_K\boldsymbol{\xi_{\sigma}}:\varepsilon(\boldsymbol{v})dx=0$ by
the inclusion property \eqref{spaces}, we have
\[
T_2=\left|\int_{\mathcal{E}_{h}}\llbracket\boldsymbol{v}\rrbracket:\{\boldsymbol{\xi_{\sigma}}\}ds\right|.
\]
Then applying the Cauchy-Schwarz  inequality yields
\[
T_2\leq\left(\int_{\mathcal{E}_{h}}C_{11}\llbracket\boldsymbol{v}\rrbracket^2ds\right)^{\frac{1}{2}}\left(\int_{\mathcal{E}_{h}}\frac{1}{C_{11}}\{\boldsymbol{\xi_{\sigma}}\}^2ds\right)^{\frac{1}{2}}
\leq \left|(\boldsymbol{\tau},\boldsymbol{v})\right|_A\left|(\boldsymbol{\xi_{\sigma}},\boldsymbol{\xi_{u}})\right|_B.
\]
Analogously, since
$\int_K\boldsymbol{\xi_{u}}\cdot(\nabla\cdot\boldsymbol{\tau})dx=0$
by \eqref{spaces}, we have by integration by parts that
\begin{align*}
T_3=\left|\int_{\mathcal{E}^i_{h}}\{\boldsymbol{\xi_{u}}\}\cdot[\boldsymbol{\tau}]ds\right|\leq\left(\int_{\mathcal{E}^i_{h}}C_{22}[\boldsymbol{\tau}]^2ds\right)^{\frac{1}{2}}\left(\int_{\mathcal{E}^i_{h}}\frac{1}{C_{22}}\{\boldsymbol{\xi_{u}}\}^2ds\right)^{\frac{1}{2}}\leq \left|(\boldsymbol{\tau},\boldsymbol{v})\right|_A\left|(\boldsymbol{\xi_{\sigma}},\boldsymbol{\xi_{u}})\right|_B.
\end{align*}
Finally, we have
\begin{align*}
T_4&=\left|\int_{\mathcal{E}_h}C_{11}\llbracket\boldsymbol{v}\rrbracket:\llbracket\boldsymbol{\xi_{u}}\rrbracket
ds\right|\leq\left(\int_{\mathcal{E}_h}C_{11}\llbracket\boldsymbol{v}\rrbracket^2ds\right)^{\frac{1}{2}}\left(\int_{\mathcal{E}_h}C_{11}\llbracket\boldsymbol{\xi_{u}}\rrbracket^2ds\right)^{\frac{1}{2}}
\leq\left|(\boldsymbol{\tau},\boldsymbol{v})\right|_A\left|(\boldsymbol{\xi_{\sigma}},\boldsymbol{\xi_{u}})\right|_B.
\end{align*}
To complete the proof, we simply have to gather the estimates of the
terms $T_i, i=1,2,3,4$, and apply once again the Cauchy-Schwarz
inequality. \quad\quad$\Box$

\begin{lem}\label{lem:3}
For any $(\boldsymbol{\sigma},\boldsymbol{u})\in
\boldsymbol{\Sigma}\times
\boldsymbol{V},(\boldsymbol{\tau},\boldsymbol{v})\in
\boldsymbol{\Sigma}_h\times \boldsymbol{V}_h$, there holds
\[
A(\boldsymbol{\tau},\boldsymbol{v};\boldsymbol{\sigma}-\boldsymbol{P}_h\boldsymbol{\sigma},\boldsymbol{u}-\boldsymbol{Q}_h\boldsymbol{u})\lesssim
\left|(\boldsymbol{\tau},\boldsymbol{v})\right|_AK_B(\boldsymbol{\sigma},\boldsymbol{u}).
\]
\end{lem}

\noindent\emph{Proof}. From the error estimates of
$\boldsymbol{P}_h$ and $\boldsymbol{Q}_h$, it follows that
\[
\left|(\boldsymbol{\sigma}-\boldsymbol{P}_h\boldsymbol{\sigma},\boldsymbol{u}-\boldsymbol{Q}_h\boldsymbol{u})\right|_B\lesssim
K_B(\boldsymbol{\sigma},\boldsymbol{u}).
\]
Together with Lemma \ref{lem:2}, we have
\begin{align*}
A(\boldsymbol{\tau},\boldsymbol{v};\boldsymbol{\sigma}-\boldsymbol{P}_h\boldsymbol{\sigma},\boldsymbol{u}-\boldsymbol{Q}_h\boldsymbol{u})&\lesssim
\left|(\boldsymbol{\tau},\boldsymbol{v})\right|_A\left|(\boldsymbol{\sigma}-\boldsymbol{P}_h\boldsymbol{\sigma},\boldsymbol{u}-\boldsymbol{Q}_h\boldsymbol{u})\right|_B
\lesssim
\left|(\boldsymbol{\tau},\boldsymbol{v})\right|_AK_B(\boldsymbol{\sigma},\boldsymbol{u}).
\quad\quad\quad\quad\quad\quad\quad\quad\quad\Box
\end{align*}

After all these preparations, we are now ready to establish an error estimate in the $A$-seminorm.
\begin{lem}\label{lem:4}
Let $(\boldsymbol{\sigma},\boldsymbol{u})\in
\boldsymbol{\Sigma}\times \boldsymbol{V}$ be the solution of
problem \eqref{eq:elas} and $(\boldsymbol{\sigma}_h,\boldsymbol{u}_h)\in
\boldsymbol{\Sigma}_h\times \boldsymbol{V}_h$ be the solution of the discrete method
\eqref{eq:variation}. Then
\begin{equation}\label{eq:Anormp}
\left|(\boldsymbol{\sigma}-\boldsymbol{\sigma}_h,\boldsymbol{u}-\boldsymbol{u}_h)\right|_A\lesssim
K^{1/2}_A(\boldsymbol{\sigma},\boldsymbol{u};\boldsymbol{\sigma},\boldsymbol{u})+K_B(\boldsymbol{\sigma},\boldsymbol{u}).
\end{equation}
\end{lem}

\noindent\emph{Proof}. It follows from the triangle inequality and Lemma \ref{lem:1} that
\begin{align}
\left|(\boldsymbol{\sigma}-\boldsymbol{\sigma}_h,\boldsymbol{u}-\boldsymbol{u}_h)\right|_A&\leq\left|(\boldsymbol{\sigma}-\boldsymbol{P}_h\boldsymbol{\sigma},\boldsymbol{u}-\boldsymbol{Q}_h\boldsymbol{u})\right|_A+\left|(\boldsymbol{P}_h\boldsymbol{\sigma}-\boldsymbol{\sigma}_h,\boldsymbol{Q}_h\boldsymbol{u}-\boldsymbol{u}_h)\right|_A
\notag \\
&\lesssim
K^{1/2}_A(\boldsymbol{\sigma},\boldsymbol{u};\boldsymbol{\sigma},\boldsymbol{u})+\left|(\boldsymbol{P}_h\boldsymbol{\sigma}-\boldsymbol{\sigma}_h,\boldsymbol{Q}_h\boldsymbol{u}-\boldsymbol{u}_h)\right|_A.
\label{eq:1}
\end{align}
By the Galerkin orthogonality \eqref{eq:go}, the definition of $A$ and
Lemma \ref{lem:3}, we see that
\begin{align*}
\left|(\boldsymbol{P}_h\boldsymbol{\sigma}-\boldsymbol{\sigma}_h,\boldsymbol{Q}_h\boldsymbol{u}-\boldsymbol{u}_h)\right|^2_A&=A(\boldsymbol{P}_h\boldsymbol{\sigma}-\boldsymbol{\sigma}_h,\boldsymbol{Q}_h\boldsymbol{u}-\boldsymbol{u}_h;\boldsymbol{P}_h\boldsymbol{\sigma}-\boldsymbol{\sigma}_h,\boldsymbol{Q}_h\boldsymbol{u}-\boldsymbol{u}_h)
\\
&=A(\boldsymbol{P}_h\boldsymbol{\sigma}-\boldsymbol{\sigma},\boldsymbol{Q}_h\boldsymbol{u}-\boldsymbol{u};\boldsymbol{P}_h\boldsymbol{\sigma}-\boldsymbol{\sigma}_h,\boldsymbol{Q}_h\boldsymbol{u}-\boldsymbol{u}_h) \\
&=A(\boldsymbol{\sigma}_h-\boldsymbol{P}_h\boldsymbol{\sigma},\boldsymbol{Q}_h\boldsymbol{u}-\boldsymbol{u}_h;\boldsymbol{\sigma}-\boldsymbol{P}_h\boldsymbol{\sigma},\boldsymbol{Q}_h\boldsymbol{u}-\boldsymbol{u})
\\
&\lesssim\left|(\boldsymbol{P}_h\boldsymbol{\sigma}-\boldsymbol{\sigma}_h,\boldsymbol{Q}_h\boldsymbol{u}-\boldsymbol{u}_h)\right|_AK_B(\boldsymbol{\sigma},\boldsymbol{u}).
\end{align*}
This implies
\begin{equation}\label{eq:3}
\left|(\boldsymbol{P}_h\boldsymbol{\sigma}-\boldsymbol{\sigma}_h,\boldsymbol{Q}_h\boldsymbol{u}-\boldsymbol{u}_h)\right|_A\lesssim
K_B(\boldsymbol{\sigma},\boldsymbol{u}).
\end{equation}
Therefore, the estimate \eqref{eq:Anormp} follows readily from \eqref{eq:1} and
\eqref{eq:3}.  \quad $\Box$ \vspace{1 pc}

We assume that the stabilization coefficients $C_{11}$ and $C_{22}$
are defined as follows:
\[
C_{11}(\boldsymbol{x})=\left\{
\begin{array}{ll}\zeta~\textrm{min}\{h_{K^+}^{\alpha_1}/p_{K^+}^{\alpha_2},h_{K^-}^{\alpha_1}/p_{K^-}^{\alpha_2}\},\quad
&\textrm{if}~ \boldsymbol{x}\in\left<K^+,K^-\right>,\\
\zeta~ h_K^{\alpha_1}/p_{K}^{\alpha_2}, &\textrm{if}~ \boldsymbol{x}\in\partial
K\cap\partial \Omega~,
\end{array}
\right.
\]
\[
C_{22}(\boldsymbol{x})=\eta~\textrm{min}\{h_{K^+}^{\beta_1}/p_{K^+}^{\beta_2},h_{K^-}^{\beta_1}/p_{K^-}^{\beta_2}\},\quad\quad
\textrm{if}~ \boldsymbol{x}\in\left<K^+,K^-\right>,
\]
with $\zeta>0,\eta>0,-1\leq\alpha_1\leq0\leq\beta_1\leq1,-1\leq\alpha_2\leq0\leq\beta_2\leq1$ independent
of the mesh size. We next introduce two symbols given by
\[
\hat{\mu}_i:=\max\{-\alpha_i,\beta_i\},\quad
\check{\mu}_i:=\min\{-\alpha_i,\beta_i\}.
\]
Denote by
\[
h:=\max\limits_{K\in\mathcal{T}_h}{h_K},\;\; p:=\min\limits_{K\in\mathcal{T}_h}{p_K},\;\; l:=\min\limits_{K\in\mathcal{T}_h}{l_K},\;\; k:=\min\limits_{K\in\mathcal{T}_h}{k_K}.
\]

\begin{thm}\label{thm:4}
Let $(\boldsymbol{\sigma},\boldsymbol{u})\in
\left(H^{1}(\Omega)\right)_{d\times d}^s\times (H^{2}(\Omega)\cap H_0^{1}(\Omega))^d$ be the solution of
problem \eqref{eq:elas} and $(\boldsymbol{\sigma}_h,\boldsymbol{u}_h)\in
\boldsymbol{\Sigma}_h\times \boldsymbol{V}_h$ be the solution of
the discrete method \eqref{eq:variation}. Assume that for each
$K\in\mathcal {T}_h$,
$(\boldsymbol{\sigma},\boldsymbol{u})|_K\in\left(H^{s+1}(K)\right)_{d\times
d}^s\times H^{s+2}(K)^d$
with integer $s\geq0$. Then for $0\leq s_{1K}\leq \min\{s,l_K\}$ and $0\leq s_{2K}\leq \min\{s+1,k_K\}$, we have
\begin{equation}\label{eq:Anorm}
\left|(\boldsymbol{\sigma}-\boldsymbol{\sigma}_h,\boldsymbol{u}-\boldsymbol{u}_h)\right|_A^2\lesssim
\sum_{K\in\mathcal{T}_h} h_K^{2\gamma_1^K}p_K^{-2\gamma_2^K}\|\boldsymbol{u}\|^2_{\max\{s_{1K}+1,s_{2K}\}+1,K},
\end{equation}
with
\[
\gamma_i^K:=\min\{s_{1K}+(1+\check{\mu}_i)/2, s_{2K}+(1-\hat{\mu}_i)/2\}, \quad i=1,2.
\]
Furthermore, for $0\leq s_{1}\leq \min\{s,l\}$ and $0\leq s_{2}\leq \min\{s+1,k\}$, there holds
\[
\left|(\boldsymbol{\sigma}-\boldsymbol{\sigma}_h,\boldsymbol{u}-\boldsymbol{u}_h)\right|_A\lesssim
h^{\gamma_1}p^{-\gamma_2} \|\boldsymbol{u}\|_{\max\{s_{1}+1,s_{2}\}+1},
\]
with
\[
\gamma_i:=\min\{s_{1}+(1+\check{\mu}_i)/2, s_{2}+(1-\hat{\mu}_i)/2\}, \quad i=1,2.
\]
\end{thm}

\noindent\emph{Proof}. From the regularity of problem
\eqref{eq:elas}, we have the regularity estimate
$\Vert\boldsymbol{\sigma}\Vert_{s+1,K}\lesssim\Vert\boldsymbol{u}\Vert_{s+2,K}$.
According to this estimate and the definition of $K_A$ we know
\begin{align}
K_A(\boldsymbol{\sigma},\boldsymbol{u};\boldsymbol{\sigma},\boldsymbol{u})&=S_1(\boldsymbol{\sigma},\boldsymbol{u};\boldsymbol{\sigma},\boldsymbol{u})+S_2(\boldsymbol{\sigma},\boldsymbol{u};\boldsymbol{\sigma},\boldsymbol{u})+S_5(\boldsymbol{\sigma},\boldsymbol{u};\boldsymbol{\sigma},\boldsymbol{u})
\notag\\
&\leq\sum_{K\in\mathcal{T}_h}\left(\frac{h_K}{p_K}\right)^{2s_{1K}+2}\|\boldsymbol{\sigma}\|^2_{s_{1K}+1,K}+
\sum_{K\in\mathcal{T}_h}C^{\partial
K}_{22}\left(\frac{h_K}{p_K}\right)^{2s_{1K}+1}\|\boldsymbol{\sigma}\|^2_{s_{1K}+1,K} \notag\\
&+\sum_{K\in\mathcal{T}_h}C^{\partial
K}_{11}\left(\frac{h_K}{p_K}\right)^{2s_{2K}+1}\|\boldsymbol{u}\|^2_{s_{2K}+1,K}\notag\\
&\leq\sum_{K\in\mathcal{T}_h}\frac{h_K^{2s_{1K}+2}}{p_K^{2s_{1K}+2}}\|\boldsymbol{u}\|^2_{s_{1K}+2,K}+
\sum_{K\in\mathcal{T}_h}\frac{h_K^{2s_{1K}+1+\beta_1}}{p_K^{2s_{1K}+1+\beta_2}}\|\boldsymbol{u}\|^2_{s_{1K}+2,K} +\sum_{K\in\mathcal{T}_h} \frac{h_K^{2s_{2K}+1+\alpha_1}}{p_K^{2s_{2K}+1+\alpha_2}}\|\boldsymbol{u}\|^2_{s_{2K}+1,K}\notag\\
&\lesssim\sum_{K\in\mathcal{T}_h} h_K^{2\gamma_1^K}p_K^{-2\gamma_2^K}\|\boldsymbol{u}\|^2_{\max\{s_{1K}+1,s_{2K}\}+1,K}.
\label{eq:6}
\end{align}
Similarly, according to the definition of $K_B$,
\begin{equation}
K^2_B(\boldsymbol{\sigma},\boldsymbol{u})\lesssim\sum_{K\in\mathcal{T}_h} h_K^{2\gamma_1^K}p_K^{-2\gamma_2^K}\|\boldsymbol{u}\|^2_{\max\{s_{1K}+1,s_{2K}\}+1,K}.
\label{eq:7}
\end{equation}
Therefore, \eqref{eq:Anorm} follows from Lemma \ref{lem:4},
\eqref{eq:6} and \eqref{eq:7}. \quad$\Box$

\begin{thm}\label{thm:5}
Suppose that $\Omega$ is a convex bounded polygon or polyhedron. Let $(\boldsymbol{\sigma},\boldsymbol{u})\in
\left(H^{1}(\Omega)\right)_{d\times d}^s\times (H^{2}(\Omega)\cap H_0^{1}(\Omega))^d$ be the solution of
problem \eqref{eq:elas} and $(\boldsymbol{\sigma}_h,\boldsymbol{u}_h)\in
\boldsymbol{\Sigma}_h\times \boldsymbol{V}_h$ be the solution of
the discrete method \eqref{eq:variation}. Assume that for each
$K\in\mathcal {T}_h$,
$(\boldsymbol{\sigma},\boldsymbol{u})|_K\in\left(H^{s+1}(K)\right)_{d\times
d}^s\times H^{s+2}(K)^d$
with integer $s\geq0$. Then for $0\leq s_{1K}\leq \min\{s,l_K\}$ and $0\leq s_{2K}\leq \min\{s+1,k_K\}$, we have
\[
\left\|\boldsymbol{u}-\boldsymbol{u}_h\right\|_0^2\lesssim
 \left(\max_{K\in\mathcal{T}_h} h_K^{2\gamma_3^K}p_K^{-2\gamma_4^K}\right)\sum_{K\in\mathcal{T}_h} h_K^{2\gamma_1^K}p_K^{-2\gamma_2^K}\|\boldsymbol{u}\|^2_{\max\{s_{1K}+1,s_{2K}\}+1,K},
\]
with
\[
\gamma_{2+i}^K:=\min\{(1+\check{\mu}_i)/2, \min\{1,k_K\}+(1-\hat{\mu}_i)/2\}, \quad i=1,2.
\]
Furthermore, for $0\leq s_{1}\leq \min\{s,l\}$ and $0\leq s_{2}\leq \min\{s+1,k\}$, there holds
\[
\left\|\boldsymbol{u}-\boldsymbol{u}_h\right\|_0\lesssim
h^{\gamma_1+\gamma_3}p^{-(\gamma_2+\gamma_4)} \|\boldsymbol{u}\|_{\max\{s_{1}+1,s_{2}\}+1},
\]
with
\[
\gamma_{2+i}:=\min\{(1+\check{\mu}_i)/2, \min\{1,k\}+(1-\hat{\mu}_i)/2\}, \quad i=1,2.
\]
\end{thm}

\noindent\emph{Proof}. We proceed by the usual duality argument. Let
$(\widetilde{\boldsymbol{\sigma}},\widetilde{\boldsymbol{u}})$ be
the solution of the auxiliary problem:
\begin{equation}\label{eq:dual1}
\left\{
\begin{array}{ll}
\mathscr{A}\boldsymbol{\widetilde{\sigma}}-\boldsymbol{\boldsymbol{\varepsilon}(\widetilde{u})}=\boldsymbol{0} & \text{in}\quad \Omega, \\
-\boldsymbol{\div} \boldsymbol{\widetilde{\sigma}}=\boldsymbol{u}-\boldsymbol{u}_h & \text{in}\quad \Omega, \\
\boldsymbol{\widetilde{u}}=\boldsymbol{0} & \text{on}\quad \partial
\Omega.
\end{array}
\right.
\end{equation}
Formally, \eqref{eq:dual1} is problem \eqref{eq:elas} with
$\boldsymbol{f}$ replaced by $\boldsymbol{u}-\boldsymbol{u}_h$. With
the same deduction as for deriving \eqref{eq:variation1}, we find
\[
A(\boldsymbol{\widetilde{\sigma}},\boldsymbol{\widetilde{u}};\boldsymbol{\tau},\boldsymbol{v})=\int_\Omega(\boldsymbol{u}-\boldsymbol{u}_h)\cdot\boldsymbol{v}dx \quad
\forall~
(\boldsymbol{\tau},\boldsymbol{v})\in\boldsymbol{\Sigma}_h\times\boldsymbol{V}_h.
\]
Now taking
$(\boldsymbol{\tau},\boldsymbol{v})=(\boldsymbol{\sigma}_h-\boldsymbol{\sigma},\boldsymbol{u}-\boldsymbol{u}_h)$,
and thanks to the definition of $A$ and the Galerkin Orthogonality
\eqref{eq:go}, we know
\begin{align}
\left\|\boldsymbol{u}-\boldsymbol{u}_h\right\|_0^2&=A(\boldsymbol{\widetilde{\sigma}},\boldsymbol{\widetilde{u}};\boldsymbol{\sigma}_h-\boldsymbol{\sigma},\boldsymbol{u}-\boldsymbol{u}_h)
=A(\boldsymbol{\sigma}-\boldsymbol{\sigma}_h,\boldsymbol{u}-\boldsymbol{u}_h;-\boldsymbol{\widetilde{\sigma}},\boldsymbol{\widetilde{u}})
\notag \\
&
=A(\boldsymbol{\sigma}-\boldsymbol{\sigma}_h,\boldsymbol{u}-\boldsymbol{u}_h;\boldsymbol{P}_h\boldsymbol{\widetilde{\sigma}}-\boldsymbol{\widetilde{\sigma}},\boldsymbol{\widetilde{u}}-\boldsymbol{Q}_h\boldsymbol{\widetilde{u}}) \notag \\
&
=A(\boldsymbol{P}_h\boldsymbol{\sigma}-\boldsymbol{\sigma}_h,\boldsymbol{Q}_h\boldsymbol{u}-\boldsymbol{u}_h;\boldsymbol{P}_h\boldsymbol{\widetilde{\sigma}}-\boldsymbol{\widetilde{\sigma}},\boldsymbol{\widetilde{u}}-\boldsymbol{Q}_h\boldsymbol{\widetilde{u}}) +A(\boldsymbol{\sigma}-\boldsymbol{P}_h\boldsymbol{\sigma},\boldsymbol{u}-\boldsymbol{Q}_h\boldsymbol{u};\boldsymbol{P}_h\boldsymbol{\widetilde{\sigma}}-\boldsymbol{\widetilde{\sigma}},\boldsymbol{\widetilde{u}}-\boldsymbol{Q}_h\boldsymbol{\widetilde{u}}).
\label{eq:2}
\end{align}
Since
$(\boldsymbol{P}_h\boldsymbol{\sigma}-\boldsymbol{\sigma}_h,\boldsymbol{Q}_h\boldsymbol{u}-\boldsymbol{u}_h)\in\boldsymbol{\Sigma}_h\times
\boldsymbol{V}_h$, it follows from Lemma \ref{lem:3} and inequality \eqref{eq:3} that
\begin{equation}\label{eq:4}
A(\boldsymbol{P}_h\boldsymbol{\sigma}-\boldsymbol{\sigma}_h,\boldsymbol{Q}_h\boldsymbol{u}-\boldsymbol{u}_h; \boldsymbol{P}_h\boldsymbol{\widetilde{\sigma}}-\boldsymbol{\widetilde{\sigma}},\boldsymbol{\widetilde{u}}-\boldsymbol{Q}_h\boldsymbol{\widetilde{u}})\lesssim
K_B(\boldsymbol{\sigma},\boldsymbol{u})K_B(-\boldsymbol{\widetilde{\sigma}},\boldsymbol{\widetilde{u}}).
\end{equation}
And according to Lemma \ref{lem:1},
\begin{equation}\label{eq:5}
A(\boldsymbol{\sigma}-\boldsymbol{P}_h\boldsymbol{\sigma},\boldsymbol{u}-\boldsymbol{Q}_h\boldsymbol{u};\boldsymbol{P}_h\boldsymbol{\widetilde{\sigma}}-\boldsymbol{\widetilde{\sigma}},\boldsymbol{\widetilde{u}}-\boldsymbol{Q}_h\boldsymbol{\widetilde{u}})\lesssim
K_A(\boldsymbol{\sigma},\boldsymbol{u};-\boldsymbol{\widetilde{\sigma}},\boldsymbol{\widetilde{u}}),
\end{equation}
which, in conjunction with \eqref{eq:2}, \eqref{eq:4} and \eqref{eq:5} implies
\begin{equation}\label{eq:8}
\left\|\boldsymbol{u}-\boldsymbol{u}_h\right\|_0^2\lesssim
K_B(\boldsymbol{\sigma},\boldsymbol{u})K_B(-\boldsymbol{\widetilde{\sigma}},\boldsymbol{\widetilde{u}})+K_A(\boldsymbol{\sigma},\boldsymbol{u};-\boldsymbol{\widetilde{\sigma}},\boldsymbol{\widetilde{u}}).
\end{equation}
By taking $s=0$ in \eqref{eq:7}, we obtain
\begin{equation}\label{eq:9}
\begin{split}
K_B(-\boldsymbol{\widetilde{\sigma}},\boldsymbol{\widetilde{u}})&\lesssim
\left(\max_{K\in\mathcal{T}_h} h_K^{\gamma_3^K}p_K^{-\gamma_4^K}\right)\|\boldsymbol{\widetilde{u}}\|_2
\lesssim \left(\max_{K\in\mathcal{T}_h} h_K^{\gamma_3^K}p_K^{-\gamma_4^K}\right)
\|\boldsymbol{u}-\boldsymbol{u}_h\|_0,
\end{split}
\end{equation}
where we have used the regularity estimate of
\eqref{eq:dual1}:~$\Vert\widetilde{\boldsymbol{u}}\Vert_2\lesssim\Vert\boldsymbol{u}-\boldsymbol{u}_h\Vert_0$.
Using the similar argument as for deriving \eqref{eq:6}, we have
\begin{align}
K_A(\boldsymbol{\sigma},\boldsymbol{u};-\boldsymbol{\widetilde{\sigma}},\boldsymbol{\widetilde{u}})\lesssim &  \left(\max_{K\in\mathcal{T}_h} h_K^{\gamma_3^K}p_K^{-\gamma_4^K}\right)\|\boldsymbol{u}-\boldsymbol{u}_h\|_0
\cdot\sqrt{\sum_{K\in\mathcal{T}_h} h_K^{2\gamma_1^K}p_K^{-2\gamma_2^K}\|\boldsymbol{u}\|^2_{\max\{s_{1K}+1,s_{2K}\}+1,K}}. \label{eq:10}
\end{align}
Finally, we can finish the proof by combining \eqref{eq:8}, \eqref{eq:9}, \eqref{eq:10} and
\eqref{eq:7} together.\quad\quad$\Box$

\begin{remark}\rm
In this paper, we derive $hp$-version error estimates for the mixed DG method \eqref{eq:mix1}-\eqref{eq:mix2} (equivalently, the method \eqref{eq:variation}) following the ideas in \cite{CastilloCockburnPerugiaSchotzau2000}. One important advantage of such arguments is that we do not require to establish the uniform inf-sup condition for the bilinear form $A$ (cf. \eqref{definition of A}). Until now, we are not able to derive such an estimate, though it plays important roles in developing a posteriori error analysis and fast solvers for this mixed DG method. It is a very challenging issue deserving further investigation.
\end{remark}

Now, let us discuss the convergence orders for some typical cases using Theorems~\ref{thm:4}-\ref{thm:5}. Write $s_{1}:=\min\{s,l\}$ and $s_{2}:=\min\{s+1,k\}$. The corresponding results are shown in Table~\ref{table:theoryfixedp}, under the condition that $(\boldsymbol{\sigma},\boldsymbol{u})\in
\left(H^{s+1}(\Omega)\right)_{d\times d}^s\times H^{s+2}(\Omega)^d$.
\begin{table}[htbp]
\tabcolsep 2pt \caption{Convergence orders in $h/p$ for some typical cases when $(\boldsymbol{\sigma},\boldsymbol{u})\in
\left(H^{s+1}(\Omega)\right)_{d\times d}^s\times H^{s+2}(\Omega)^d$,
$s\geq 0.$}\label{table:theoryfixedp} \vspace*{-8pt}
\begin{center}
\def\temptablewidth{12.6cm}
\begin{tabular*}{\temptablewidth}{@{\extracolsep{\fill}}|c|c|c|c|}
\hline $C_{11}$ & $C_{22}$ &
$\left|(\boldsymbol{\sigma}-\boldsymbol{\sigma}_h,\boldsymbol{u}-\boldsymbol{u}_h)\right|_A$
& $\left\|\boldsymbol{u}-\boldsymbol{u}_h\right\|_0$ \\
\hline $O(p/h)$ & $O(1)$ & $\min\{s+1/2,l+1/2,k\}$&
$\min\{s+1/2,l+1/2,k\}+\min\{1/2,k\}$
\\
\hline $O(p/h)$ & $O(h/p)$ & $\min\{s+1,l+1,k\}$ &
$\min\{s+1,l+1,k\}+\min\{1,k\}$
\\
\hline $O(1)$ & $O(1)$ & $\min\{s,l,k\}+1/2$& $\min\{s,l,k\}+1$
\\
\hline $O(1)$ & $O(h/p)$ & $\min\{s+1/2,l+1/2,k\}$&
$\min\{s+1/2,l+1/2,k\}+\min\{1/2,k\}$
\\
\hline
       \end{tabular*}
       \end{center}
       \end{table}

\begin{remark}\rm
According to the $h$-version error estimates in \cite{ChenHuangHuangXu2010} and Table~\ref{table:theoryfixedp}, we can find that
the convergence rates of the errors in $L^2$ norm and energy norm with $C_{11}=O(h^{-1})$ and $C_{22}=0$ (i.e., the LDG method in \cite{ChenHuangHuangXu2010}) coincide with those of the DG method \eqref{eq:mix1}-\eqref{eq:mix2} with $C_{11}=O(h^{-1})$ and $C_{22}=O(h)$. We will also observe this phenomenon from the numerical experiments in Section 4.
\end{remark}

\section{Numerical results}

In this section, we intend to present a variety of numerical examples in order to illustrate the numerical performance of the mixed DG method \eqref{eq:variation} (or equivalently, the method \eqref{eq:mix1}-\eqref{eq:mix2}). In all the numerical examples, we choose $\lambda=0.3$ and $\mu=0.35$. For
any $K\in \mathcal{T}_h$, we take $\mathcal{S}_1(K)=P_l(K)$ and $\mathcal
{S}_2(K)=P_k(K)$ where $k,l\geq0$. Set $\eta=1$ when $C_{22}\neq0$, and let $\zeta=1$.

\subsection{A two-dimensional example}
Let
$\Omega=(-1,1)\times(-1,1)$,  and
\begin{align*}
\boldsymbol{f}(x_1,x_2)=&\left(\begin{array}{l}
-8(x_1 + x_2)\left((3x_1x_2-2)(x_1^2+x_2^2)+5(x_1x_2-1)^2-2x_1^2x_2^2\right) \\
-8(x_1 - x_2)\left((3x_1x_2+2)(x_1^2+x_2^2)-5(x_1x_2+1)^2+2x_1^2x_2^2\right)
\end{array}\right).
\end{align*}
It can be verified the exact solution of \eqref{eq:elas} is
\[
\boldsymbol{u}(x_1,x_2)=\frac{80}{7}\left(\begin{array}{l}
-x_2(1-x_2^2)(1-x_1^2)^2 \\
x_1(1-x_1^2)(1-x_2^2)^2
\end{array}\right)-4\left(\begin{array}{l}
x_1(1-x_1^2)(1-x_2^2)^2 \\
x_2(1-x_2^2)(1-x_1^2)^2
\end{array}\right).
\]

First of all,  we use the uniform triangulation $\mathcal{T}_h$ of $\Omega$ and consider the $h$-version convergence of our DG method with fixed $p$. In this example, we take $\alpha_2=\beta_2=0$.
Tables~\ref{table:errorL2k0l0}-\ref{table:errorEnergyk0l0} show the errors in $L^2$ norm and energy norm for $k=0, l=0$, respectively. It is observed from Tables~\ref{table:errorL2k0l0}-\ref{table:errorEnergyk0l0} that the numerical convergence rates of $\|\boldsymbol{u}-\boldsymbol{u}_h\|_0$ and $\left|(\boldsymbol{\sigma}-\boldsymbol{\sigma}_h,\boldsymbol{u}-\boldsymbol{u}_h)\right|_A$ agree with the theoretical convergence rates in Table~\ref{table:theoryfixedp} except the case $C_{11}=O(1), C_{22}=O(h)$, in which both the numerical convergence rates of $\|\boldsymbol{u}-\boldsymbol{u}_h\|_0$ and $\left|(\boldsymbol{\sigma}-\boldsymbol{\sigma}_h,\boldsymbol{u}-\boldsymbol{u}_h)\right|_A$ are half-order higher than the theoretical results.
And numerical results for $k=1, l=0$ given in Tables~\ref{table:errorL2k1l0}-\ref{table:errorEnergyk1l0} illustrate that numerical convergence rates of $\|\boldsymbol{u}-\boldsymbol{u}_h\|_0$ and $\left|(\boldsymbol{\sigma}-\boldsymbol{\sigma}_h,\boldsymbol{u}-\boldsymbol{u}_h)\right|_A$ are all consistent with the theoretical results.
Tables~\ref{table:errorL2k1l1}-\ref{table:errorEnergyk1l1} present the errors in $L^2$ norm and energy norm for $k=1, l=1$, respectively,
from which we can see that $\|\boldsymbol{u}-\boldsymbol{u}_h\|_0=O(h^2)$ for different choices of $C_{11}$ and $C_{22}$, and $\left|(\boldsymbol{\sigma}-\boldsymbol{\sigma}_h,\boldsymbol{u}-\boldsymbol{u}_h)\right|_A=O(h)$ for $C_{22}=O(h), 0$ and $\left|(\boldsymbol{\sigma}-\boldsymbol{\sigma}_h,\boldsymbol{u}-\boldsymbol{u}_h)\right|_A=O(h^{3/2})$ for $C_{22}=O(1)$. Thus,
the numerical convergence rates of $\|\boldsymbol{u}-\boldsymbol{u}_h\|_0$  coincide with the theoretical results in Theorem~\ref{thm:5} except the cases $C_{11}=O(h^{-1}), C_{22}=O(1)$ and $C_{11}=O(1), C_{22}=O(h)$, in which the numerical convergence rates of $\|\boldsymbol{u}-\boldsymbol{u}_h\|_0$ are half-order higher than the theoretical convergence rates. The numerical convergence rates of $\left|(\boldsymbol{\sigma}-\boldsymbol{\sigma}_h,\boldsymbol{u}-\boldsymbol{u}_h)\right|_A$  also coincide with the theoretical results in Theorem~\ref{thm:4} except the case $C_{11}=O(h^{-1}), C_{22}=O(1)$, in which the numerical convergence rate of $\left|(\boldsymbol{\sigma}-\boldsymbol{\sigma}_h,\boldsymbol{u}-\boldsymbol{u}_h)\right|_A$ is half-order higher than the theoretical convergence rate.
The numerical errors $\|\boldsymbol{u}-\boldsymbol{u}_h\|_0$ and $\left|(\boldsymbol{\sigma}-\boldsymbol{\sigma}_h,\boldsymbol{u}-\boldsymbol{u}_h)\right|_A$ for $k=2, l=2$ are listed in Tables~\ref{table:errorL2k2l2}-\ref{table:errorEnergyk2l2}, respectively.
We can see that $\|\boldsymbol{u}-\boldsymbol{u}_h\|_0=O(h^3)$ for different choices of $C_{11}$ and $C_{22}$, and $\left|(\boldsymbol{\sigma}-\boldsymbol{\sigma}_h,\boldsymbol{u}-\boldsymbol{u}_h)\right|_A=O(h^{5/2})$ for $C_{11}=O(1), C_{22}=O(1)$ and
$\left|(\boldsymbol{\sigma}-\boldsymbol{\sigma}_h,\boldsymbol{u}-\boldsymbol{u}_h)\right|_A=O(h^{2})$ for other choices of $C_{11}$ and $C_{22}$. Again,
the numerical convergence rates of $\|\boldsymbol{u}-\boldsymbol{u}_h\|_0$  coincide with the theoretical results in Theorem~\ref{thm:5} except the cases $C_{11}=O(h^{-1}), C_{22}=O(1)$ and $C_{11}=O(1), C_{22}=O(h)$, whose numerical convergence rates are half-order higher than the theoretical convergence rates. Whereas all the numerical convergence rates of $\left|(\boldsymbol{\sigma}-\boldsymbol{\sigma}_h,\boldsymbol{u}-\boldsymbol{u}_h)\right|_A$ coincide with the theoretical results in Theorem~\ref{thm:4} for the different choices of $C_{11}$ and $C_{22}$.
We also list the numerical errors $\|\boldsymbol{u}-\boldsymbol{u}_h\|_0$ and $\left|(\boldsymbol{\sigma}-\boldsymbol{\sigma}_h,\boldsymbol{u}-\boldsymbol{u}_h)\right|_A$  for $k=2, l=2$, $C_{11}=O(h^{-1}), C_{22}=O(h^{-1})$ in Table~\ref{table:errorsk2l2C11C22hm1}. The choice $C_{11}=O(h^{-1})$ and $C_{22}=O(h^{-1})$ is not covered by our theoretical analysis. Both the numerical errors $\|\boldsymbol{u}-\boldsymbol{u}_h\|_0$ and $\left|(\boldsymbol{\sigma}-\boldsymbol{\sigma}_h,\boldsymbol{u}-\boldsymbol{u}_h)\right|_A$ approximate $O(h^2)$, which implies the numerical convergence rate of $\|\boldsymbol{u}-\boldsymbol{u}_h\|_0$ is one-order lower than those for the choices of $C_{11}$ and $C_{22}$ in Table~\ref{table:errorL2k2l2}. 
\begin{table}[htbp]
\centering
\tabcolsep 2pt \caption{Uniform triangular meshes: Error $\|\boldsymbol{u}-\boldsymbol{u}_h\|_0$ vs $h$
for different choices of $C_{11}, C_{22}$ when
$k=0, l=0$.}\label{table:errorL2k0l0} \vspace*{3pt}
\resizebox{0.98\textwidth}{!}{
\def\temptablewidth{\textwidth}
\begin{tabular*}{\temptablewidth}{@{\extracolsep{\fill}}|c|c|c|c|c|c|c|c|c|c|c|}
\hline \backslashbox{\kern3em$h$\kern-3em}{\kern2em$C_{11}, C_{22}$} & $O(h^{-1}), 0$ & order & $O(h^{-1}), O(1)$ & order & $O(h^{-1}), O(h)$ & order & $O(1), O(1)$ & order & $O(1), O(h)$ & order \\
\hline $1$ & 4.1084E+00 & $-$ & 5.8643E+00 & $-$ & 6.0867E+00 & $-$ & 5.2390E+00 & $-$ & 5.4575E+00 & $-$ \\
\hline $2^{-1}$ & 3.1212E+00 & 0.40 & 3.1782E+00 & 0.88 & 3.0151E+00 & 1.01 & 4.5372E+00 & 0.21 & 3.7891E+00 & 0.53 \\
\hline $2^{-2}$ & 2.9507E+00 & 0.08 & 2.4644E+00 & 0.37 & 2.6930E+00 & 0.16 & 3.0346E+00 & 0.58 & 1.9067E+00 & 0.99 \\
\hline $2^{-3}$ & 2.9796E+00 & 0 & 2.6170E+00 & 0 & 2.8852E+00 & 0 & 1.7554E+00 & 0.79 & 9.1402E$-$01 & 1.06 \\
\hline $2^{-4}$ & 3.0255E+00 & 0 & 2.8213E+00 & 0 & 2.9964E+00 & 0 & 9.4401E$-$01 & 0.90 & 4.7760E$-$01 & 0.94 \\
\hline $2^{-5}$ & 3.0553E+00 & 0 & 2.9481E+00 & 0 & 3.0460E+00 & 0 & 4.9068E$-$01 & 0.94 & 2.5513E$-$01 & 0.91 \\
\hline $2^{-6}$ & 3.0717E+00 & 0 & 3.0169E+00 & 0 & 3.0684E+00 & 0 & 2.5064E$-$01 & 0.97 & 1.3345E$-$01 & 0.94 \\
\hline
       \end{tabular*}
       }
\end{table}
\begin{table}[htbp]
\centering
\tabcolsep 2pt \caption{Uniform triangular meshes: Error $\left|(\boldsymbol{\sigma}-\boldsymbol{\sigma}_h,\boldsymbol{u}-\boldsymbol{u}_h)\right|_A$  vs $h$
for different choices of $C_{11}, C_{22}$ when
$k=0, l=0$.}\label{table:errorEnergyk0l0} \vspace*{3pt}
\resizebox{0.98\textwidth}{!}{
\def\temptablewidth{\textwidth}
\begin{tabular*}{\temptablewidth}{@{\extracolsep{\fill}}|c|c|c|c|c|c|c|c|c|c|c|}
\hline \backslashbox{\kern3em$h$\kern-3em}{\kern2em$C_{11}, C_{22}$} & $O(h^{-1}), 0$ & order & $O(h^{-1}), O(1)$ & order & $O(h^{-1}), O(h)$ & order & $O(1), O(1)$ & order & $O(1), O(h)$ & order \\
\hline $1$ & 1.4858E+01 & $-$ & 1.9821E+01 & $-$ & 2.0202E+01 & $-$ & 1.9383E+01 & $-$ & 1.9761E+01 & $-$ \\
\hline $2^{-1}$ & 1.1910E+01 & 0.32 & 1.6989E+01 & 0.22 & 1.5638E+01 & 0.37 & 1.8133E+01 & 0.10 & 1.6306E+01 & 0.28 \\
\hline $2^{-2}$ & 1.0314E+01 & 0.21 & 1.4139E+01 & 0.27 & 1.1905E+01 & 0.39 & 1.4434E+01 & 0.33 & 1.0633E+01 & 0.62 \\
\hline $2^{-3}$ & 9.8875E+00 & 0.06 & 1.2325E+01 & 0.20 & 1.0373E+01 & 0.20 & 1.0780E+01 & 0.42 & 6.5685E+00 & 0.70 \\
\hline $2^{-4}$ & 9.8475E+00 & 0.01 & 1.1261E+01 & 0.13 & 9.9699E+00 & 0.06 & 7.8341E+00 & 0.46 & 4.1821E+00 & 0.65 \\
\hline $2^{-5}$ & 9.8769E+00 & 0 & 1.0655E+01 & 0.08 & 9.9039E+00 & 0.01 & 5.6192E+00 & 0.48 & 2.7813E+00 & 0.59 \\
\hline $2^{-6}$ & 9.9041E+00 & 0 & 1.0319E+01 & 0.05 & 9.9087E+00 & 0 & 4.0031E+00 & 0.49 & 1.9054E+00 & 0.55 \\
\hline
       \end{tabular*}
       }
\end{table}
\begin{table}[htbp]
\centering
\tabcolsep 2pt \caption{Uniform triangular meshes: Error $\|\boldsymbol{u}-\boldsymbol{u}_h\|_0$ vs $h$
for different choices of $C_{11}, C_{22}$ when
$k=1, l=0$.}\label{table:errorL2k1l0} \vspace*{3pt}
\resizebox{0.98\textwidth}{!}{
\def\temptablewidth{\textwidth}
\begin{tabular*}{\temptablewidth}{@{\extracolsep{\fill}}|c|c|c|c|c|c|c|c|c|c|c|}
\hline \backslashbox{\kern3em$h$\kern-3em}{\kern2em$C_{11}, C_{22}$} & $O(h^{-1}), 0$ & order & $O(h^{-1}), O(1)$ & order & $O(h^{-1}), O(h)$ & order & $O(1), O(1)$ & order & $O(1), O(h)$ & order \\
\hline $1$ & 6.3925E+00 & $-$ & 1.0176E+01 & $-$ & 1.0635E+01 & $-$ & 9.9955E+00 & $-$ & 1.0468E+01 & $-$ \\
\hline $2^{-1}$ & 2.1140E+00 & 1.60 & 6.2118E+00 & 0.71 & 4.0434E+00 & 1.40 & 7.5757E+00 & 0.40 & 5.6379E+00 & 0.89 \\
\hline $2^{-2}$ & 4.9023E$-$01 & 2.11 & 3.6613E+00 & 0.76 & 1.1967E+00 & 1.76 & 4.3381E+00 & 0.80 & 2.2842E+00 & 1.30 \\
\hline $2^{-3}$ & 1.0629E$-$01 & 2.21 & 2.0404E+00 & 0.84 & 3.2071E$-$01 & 1.90 & 2.3415E+00 & 0.89 & 8.9488E$-$01 & 1.35 \\
\hline $2^{-4}$ & 2.3212E$-$02 & 2.20 & 1.1004E+00 & 0.89 & 8.2405E$-$02 & 1.96 & 1.2363E+00 & 0.92 & 3.5130E$-$01 & 1.35 \\
\hline $2^{-5}$ & 5.2841E$-$03 & 2.14 & 5.7790E$-$01 & 0.93 & 2.0823E$-$02 & 1.98 & 6.4134E$-$01 & 0.95 & 1.4024E$-$01 & 1.32 \\
\hline $2^{-6}$ & 1.2504E$-$03 & 2.08 & 2.9754E$-$01 & 0.96 & 5.2292E$-$03 & 1.99 & 3.2802E$-$01 & 0.97 & 5.8102E$-$02 & 1.27 \\
\hline $2^{-7}$ & 3.0347E$-$04 & 2.04 & 1.5122E$-$01 & 0.98 & 1.3099E$-$03 & 2.00 & 1.6613E$-$01 & 0.98 & 2.5314E$-$02 & 1.20 \\
\hline $2^{-8}$ & 7.4710E$-$05 & 2.02 & 7.6271E$-$02 & 0.99 & 3.2779E$-$04 & 2.00 & 8.3640E$-$02 & 0.99 & 1.1580E$-$02 & 1.13 \\
\hline
       \end{tabular*}
       }
\end{table}
\begin{table}[htbp]
\centering
\tabcolsep 2pt \caption{Uniform triangular meshes: Error $\left|(\boldsymbol{\sigma}-\boldsymbol{\sigma}_h,\boldsymbol{u}-\boldsymbol{u}_h)\right|_A$  vs $h$
for different choices of $C_{11}, C_{22}$ when
$k=1, l=0$.}\label{table:errorEnergyk1l0} \vspace*{3pt}
\resizebox{0.98\textwidth}{!}{
\def\temptablewidth{\textwidth}
\begin{tabular*}{\temptablewidth}{@{\extracolsep{\fill}}|c|c|c|c|c|c|c|c|c|c|c|}
\hline \backslashbox{\kern3em$h$\kern-3em}{\kern2em$C_{11}, C_{22}$} & $O(h^{-1}), 0$ & order & $O(h^{-1}), O(1)$ & order & $O(h^{-1}), O(h)$ & order & $O(1), O(1)$ & order & $O(1), O(h)$ & order \\
\hline $1$ & 1.9375E+01 & $-$ & 2.5339E+01 & $-$ & 2.5783E+01 & $-$ & 2.5230E+01 & $-$ & 2.5681E+01 & $-$ \\
\hline $2^{-1}$ & 1.1418E+01 & 0.76 & 1.8456E+01 & 0.46 & 1.5920E+01 & 0.70 & 1.9899E+01 & 0.34 & 1.7428E+01 & 0.56 \\
\hline $2^{-2}$ & 5.6459E+00 & 1.02 & 1.3170E+01 & 0.49 & 8.5682E+00 & 0.89 & 1.4264E+01 & 0.48 & 1.0055E+01 & 0.79 \\
\hline $2^{-3}$ & 2.7144E+00 & 1.06 & 9.5205E+00 & 0.47 & 4.4004E+00 & 0.96 & 1.0175E+01 & 0.49 & 5.6226E+00 & 0.84 \\
\hline $2^{-4}$ & 1.3113E+00 & 1.05 & 6.8917E+00 & 0.47 & 2.2215E+00 & 0.99 & 7.2687E+00 & 0.49 & 3.1670E+00 & 0.83 \\
\hline $2^{-5}$ & 6.4119E$-$01 & 1.03 & 4.9574E+00 & 0.48 & 1.1147E+00 & 1.00 & 5.1831E+00 & 0.49 & 1.8415E+00 & 0.78 \\
\hline $2^{-6}$ & 3.1657E$-$01 & 1.02 & 3.5426E+00 & 0.49 & 5.5814E$-$01 & 1.00 & 3.6854E+00 & 0.49 & 1.1235E+00 & 0.71 \\
\hline $2^{-7}$ & 1.5723E$-$01 & 1.01 & 2.5198E+00 & 0.49 & 2.7924E$-$01 & 1.00 & 2.6145E+00 & 0.50 & 7.2037E$-$01 & 0.64 \\
\hline $2^{-8}$ & 7.8340E$-$02 & 1.01 & 1.7874E+00 & 0.50 & 1.3966E$-$01 & 1.00 & 1.8520E+00 & 0.50 & 4.8056E$-$01 & 0.58 \\
\hline
       \end{tabular*}
       }
\end{table}
\begin{table}[htbp]
\centering
\tabcolsep 2pt
\caption{Uniform triangular meshes: Error $\|\boldsymbol{u}-\boldsymbol{u}_h\|_0$ vs $h$
for different choices of $C_{11}, C_{22}$ when
$k=1, l=1$.}\label{table:errorL2k1l1} \vspace*{3pt}
\resizebox{0.98\textwidth}{!}{
\def\temptablewidth{\textwidth}
\begin{tabular*}{\temptablewidth}{@{\extracolsep{\fill}}|c|c|c|c|c|c|c|c|c|c|c|}
\hline \backslashbox{\kern3em$h$\kern-3em}{\kern2em$C_{11}, C_{22}$} & $O(h^{-1}), 0$ & order & $O(h^{-1}), O(1)$ & order & $O(h^{-1}), O(h)$ & order & $O(1), O(1)$ & order & $O(1), O(h)$ & order \\
\hline $1$ & 2.7907E+00 & $-$ & 3.7355E+00 & $-$ & 3.9650E+00 & $-$ & 3.7263E+00 & $-$ & 3.9531E+00 & $-$ \\
\hline $2^{-1}$ & 1.0974E+00 & 1.35 & 9.8826E$-$01 & 1.92 & 9.2265E$-$01 & 2.10 & 1.1226E+00 & 1.73 & 1.0130E+00 & 1.96 \\
\hline $2^{-2}$ & 3.8060E$-$01 & 1.53 & 2.4488E$-$01 & 2.01 & 2.2904E$-$01 & 2.01 & 2.8496E$-$01 & 1.98 & 2.4339E$-$01 & 2.06 \\
\hline $2^{-3}$ & 1.1402E$-$01 & 1.74 & 6.0548E$-$02 & 2.02 & 5.7465E$-$02 & 1.99 & 7.0305E$-$02 & 2.02 & 5.7844E$-$02 & 2.07 \\
\hline $2^{-4}$ & 3.0809E$-$02 & 1.89 & 1.5001E$-$02 & 2.01 & 1.4395E$-$02 & 2.00 & 1.7475E$-$02 & 2.01 & 1.3818E$-$02 & 2.07 \\
\hline $2^{-5}$ & 7.9521E$-$03 & 1.95 & 3.7250E$-$03 & 2.01 & 3.6014E$-$03 & 2.00 & 4.3787E$-$03 & 2.00 & 3.3506E$-$03 & 2.04 \\
\hline $2^{-6}$ & 2.0155E$-$03 & 1.98 & 9.2743E$-$04 & 2.01 & 9.0061E$-$04 & 2.00 & 1.1009E$-$03 & 1.99 & 8.2381E$-$04 & 2.02 \\
\hline $2^{-7}$ & 5.0703E$-$04 & 1.99 & 2.3134E$-$04 & 2.00 & 2.2517E$-$04 & 2.00 & 2.7694E$-$04 & 1.99 & 2.0452E$-$04 & 2.01 \\
\hline
       \end{tabular*}
       }
\end{table}
\begin{table}[htbp]
\centering
\tabcolsep 2pt \caption{Uniform triangular meshes: Error $\left|(\boldsymbol{\sigma}-\boldsymbol{\sigma}_h,\boldsymbol{u}-\boldsymbol{u}_h)\right|_A$  vs $h$
for different choices of $C_{11}, C_{22}$ when
$k=1, l=1$.}\label{table:errorEnergyk1l1} \vspace*{3pt}
\resizebox{0.98\textwidth}{!}{
\def\temptablewidth{\textwidth}
\begin{tabular*}{\temptablewidth}{@{\extracolsep{\fill}}|c|c|c|c|c|c|c|c|c|c|c|}
\hline \backslashbox{\kern3em$h$\kern-3em}{\kern2em$C_{11}, C_{22}$} & $O(h^{-1}), 0$ & order & $O(h^{-1}), O(1)$ & order & $O(h^{-1}), O(h)$ & order & $O(1), O(1)$ & order & $O(1), O(h)$ & order \\
\hline $1$ & 9.0987E+00 & $-$ & 1.2493E+01 & $-$ & 1.2836E+01 & $-$ & 1.2471E+01 & $-$ & 1.2812E+01 & $-$ \\
\hline $2^{-1}$ & 5.3420E+00 & 0.77 & 4.8502E+00 & 1.37 & 4.8539E+00 & 1.40 & 4.9027E+00 & 1.35 & 4.7674E+00 & 1.43 \\
\hline $2^{-2}$ & 3.2250E+00 & 0.73 & 1.6528E+00 & 1.55 & 2.0256E+00 & 1.26 & 1.7202E+00 & 1.51 & 1.8064E+00 & 1.40 \\
\hline $2^{-3}$ & 1.7963E+00 & 0.84 & 5.4964E$-$01 & 1.59 & 9.3573E$-$01 & 1.11 & 5.9770E$-$01 & 1.53 & 7.4950E$-$01 & 1.27 \\
\hline $2^{-4}$ & 9.4376E$-$01 & 0.93 & 1.8558E$-$01 & 1.57 & 4.5545E$-$01 & 1.04 & 2.1027E$-$01 & 1.51 & 3.3427E$-$01 & 1.16 \\
\hline $2^{-5}$ & 4.8247E$-$01 & 0.97 & 6.3630E$-$02 & 1.54 & 2.2579E$-$01 & 1.01 & 7.4592E$-$02 & 1.50 & 1.5596E$-$01 & 1.10 \\
\hline $2^{-6}$ & 2.4373E$-$01 & 0.99 & 2.2080E$-$02 & 1.53 & 1.1257E$-$01 & 1.00 & 2.6544E$-$02 & 1.49 & 7.4944E$-$02 & 1.06 \\
\hline $2^{-7}$ & 1.2246E$-$01 & 0.99 & 7.7249E$-$03 & 1.52 & 5.6229E$-$02 & 1.00 & 9.4483E$-$03 & 1.49 & 3.6683E$-$02 & 1.03 \\
\hline
       \end{tabular*}
       }
\end{table}
\begin{table}[htbp]
\centering
\tabcolsep 2pt
\caption{Uniform triangular meshes: Error $\|\boldsymbol{u}-\boldsymbol{u}_h\|_0$ vs $h$
for different choices of $C_{11}, C_{22}$ when
$k=2, l=2$.}\label{table:errorL2k2l2} \vspace*{3pt}
\resizebox{0.98\textwidth}{!}{
\def\temptablewidth{\textwidth}
\begin{tabular*}{\temptablewidth}{@{\extracolsep{\fill}}|c|c|c|c|c|c|c|c|c|c|c|}
\hline \backslashbox{\kern3em$h$\kern-3em}{\kern2em$C_{11}, C_{22}$} & $O(h^{-1}), 0$ & order & $O(h^{-1}), O(1)$ & order & $O(h^{-1}), O(h)$ & order & $O(1), O(1)$ & order & $O(1), O(h)$ & order \\
\hline $1$ & 1.1360E+00 & $-$ & 1.0723E+00 & $-$ & 1.0858E+00 & $-$ & 1.0610E+00 & $-$ & 1.0742E+00 & $-$ \\
\hline $2^{-1}$ & 1.5057E$-$01 & 2.92 & 1.4285E$-$01 & 2.91 & 1.3748E$-$01 & 2.98 & 1.5784E$-$01 & 2.75 & 1.4760E$-$01 & 2.86 \\
\hline $2^{-2}$ & 1.8804E$-$02 & 3.00 & 1.8380E$-$02 & 2.96 & 1.7271E$-$02 & 2.99 & 2.2030E$-$02 & 2.84 & 1.8183E$-$02 & 3.02 \\
\hline $2^{-3}$ & 2.3375E$-$03 & 3.01 & 2.3432E$-$03 & 2.97 & 2.1419E$-$03 & 3.01 & 2.9393E$-$03 & 2.91 & 2.1962E$-$03 & 3.05 \\
\hline $2^{-4}$ & 2.9170E$-$04 & 3.00 & 2.9842E$-$04 & 2.97 & 2.6640E$-$04 & 3.01 & 3.8352E$-$04 & 2.94 & 2.7242E$-$04 & 3.01 \\
\hline $2^{-5}$ & 3.6456E$-$05 & 3.00 & 3.7896E$-$05 & 2.98 & 3.3229E$-$05 & 3.00 & 4.9268E$-$05 & 2.96 & 3.4244E$-$05 & 2.99 \\
\hline $2^{-6}$ & 4.5574E$-$06 & 3.00 & 4.7977E$-$06 & 2.98 & 4.1500E$-$06 & 3.00 & 6.2654E$-$06 & 2.98 & 4.3088E$-$06 & 2.99 \\
\hline
       \end{tabular*}
       }
\end{table}
\begin{table}[htbp]
\centering
\tabcolsep 2pt \caption{Uniform triangular meshes: Error $\left|(\boldsymbol{\sigma}-\boldsymbol{\sigma}_h,\boldsymbol{u}-\boldsymbol{u}_h)\right|_A$  vs $h$
for different choices of $C_{11}, C_{22}$ when
$k=2, l=2$.}\label{table:errorEnergyk2l2} \vspace*{3pt}
\resizebox{0.98\textwidth}{!}{
\def\temptablewidth{\textwidth}
\begin{tabular*}{\temptablewidth}{@{\extracolsep{\fill}}|c|c|c|c|c|c|c|c|c|c|c|}
\hline \backslashbox{\kern3em$h$\kern-3em}{\kern2em$C_{11}, C_{22}$} & $O(h^{-1}), 0$ & order & $O(h^{-1}), O(1)$ & order & $O(h^{-1}), O(h)$ & order & $O(1), O(1)$ & order & $O(1), O(h)$ & order \\
\hline $1$ & 4.1805E+00 & $-$ & 3.9577E+00 & $-$ & 3.9892E+00 & $-$ & 3.9340E+00 & $-$ & 3.9650E+00 & $-$ \\
\hline $2^{-1}$ & 1.2298E+00 & 1.77 & 9.7099E$-$01 & 2.03 & 9.8231E$-$01 & 2.02 & 9.3959E$-$01 & 2.07 & 9.2361E$-$01 & 2.10 \\
\hline $2^{-2}$ & 3.1918E$-$01 & 1.95 & 2.2150E$-$01 & 2.13 & 2.3515E$-$01 & 2.06 & 1.8159E$-$01 & 2.37 & 1.7344E$-$01 & 2.41 \\
\hline $2^{-3}$ & 8.1303E$-$02 & 1.97 & 5.3238E$-$02 & 2.06 & 5.8599E$-$02 & 2.00 & 3.2842E$-$02 & 2.47 & 3.2735E$-$02 & 2.41 \\
\hline $2^{-4}$ & 2.0655E$-$02 & 1.98 & 1.3188E$-$02 & 2.01 & 1.4725E$-$02 & 1.99 & 5.8450E$-$03 & 2.49 & 6.7115E$-$03 & 2.29 \\
\hline $2^{-5}$ & 5.2190E$-$03 & 1.98 & 3.2964E$-$03 & 2.00 & 3.6962E$-$03 & 1.99 & 1.0370E$-$03 & 2.49 & 1.4869E$-$03 & 2.17 \\
\hline $2^{-6}$ & 1.3127E$-$03 & 1.99 & 8.2533E$-$04 & 2.00 & 9.2623E$-$04 & 2.00 & 1.8383E$-$04 & 2.50 & 3.4748E$-$04 & 2.10 \\
\hline
       \end{tabular*}
       }
\end{table}
\begin{table}[htbp]
\centering
\tabcolsep 2pt \caption{Uniform triangular meshes: Errors $\|\boldsymbol{u}-\boldsymbol{u}_h\|_0$ and $\left|(\boldsymbol{\sigma}-\boldsymbol{\sigma}_h,\boldsymbol{u}-\boldsymbol{u}_h)\right|_A$  vs $h$
when $k=2, l=2$, $C_{11}=O(h^{-1})$ and $C_{22}=O(h^{-1})$.}\label{table:errorsk2l2C11C22hm1} \vspace*{3pt}
\def\temptablewidth{0.45\textwidth}
\begin{tabular*}{\temptablewidth}{@{\extracolsep{\fill}}|c|c|c|c|c|}
\hline $h$ & $\|\boldsymbol{u}-\boldsymbol{u}_h\|_0$ & order & $\left|(\boldsymbol{\sigma}-\boldsymbol{\sigma}_h,\boldsymbol{u}-\boldsymbol{u}_h)\right|_A$ & order \\
\hline $1$ & 1.0592E+00 & $-$ & 3.9225E+00 & $-$ \\
\hline $2^{-1}$ & 1.4998E$-$01 & 2.82 & 9.6951E$-$01 & 2.02  \\
\hline $2^{-2}$ & 2.0049E$-$02 & 2.90 & 2.1581E$-$01 & 2.17  \\
\hline $2^{-3}$ & 3.0847E$-$03 & 2.70 & 5.1020E$-$02 & 2.08  \\
\hline $2^{-4}$ & 6.0357E$-$04 & 2.35 & 1.2557E$-$02 & 2.02  \\
\hline $2^{-5}$ & 1.3910E$-$04 & 2.12 & 3.1295E$-$03 & 2.00  \\
\hline
       \end{tabular*}
\end{table}

\begin{figure}
  \centering
  \includegraphics[width=5cm]{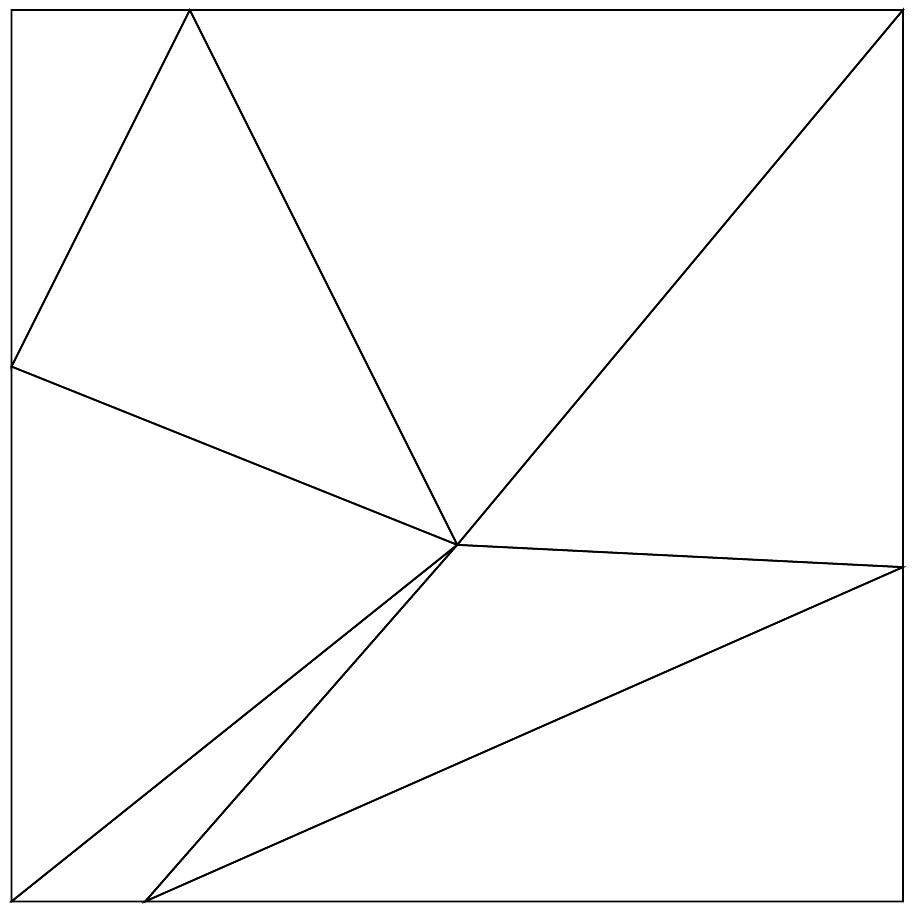}\\
  \caption{The initial unstructured triangular mesh.}\label{fig:initialunstructuredmesh}
\end{figure}

Then we subdivide the domain $\Omega$ by the unstructured triangular meshes. The initial mesh is shown in Fig.~\ref{fig:initialunstructuredmesh}, which is further refined by connecting midpoints of all edges of each triangle.
The numerical errors $\|\boldsymbol{u}-\boldsymbol{u}_h\|_0$ and $\left|(\boldsymbol{\sigma}-\boldsymbol{\sigma}_h,\boldsymbol{u}-\boldsymbol{u}_h)\right|_A$ for $k=1, l=1$ are presented in Tables~\ref{table:errorL2k1l1unstructure}-\ref{table:errorEnergyk1l1unstructure}, respectively.
Compared with Tables~\ref{table:errorL2k1l1}-\ref{table:errorEnergyk1l1}, the numerical convergence rates indicated in Tables~\ref{table:errorL2k1l1unstructure}-\ref{table:errorEnergyk1l1unstructure} are the same as those on the uniform triangular meshes.
That is,
the numerical convergence rates of $\|\boldsymbol{u}-\boldsymbol{u}_h\|_0$ coincide with the theoretical results except the cases $C_{11}=O(h^{-1}), C_{22}=O(1)$ and $C_{11}=O(1), C_{22}=O(h)$. And the numerical convergence rates of $\left|(\boldsymbol{\sigma}-\boldsymbol{\sigma}_h,\boldsymbol{u}-\boldsymbol{u}_h)\right|_A$ coincide with the theoretical results except the case $C_{11}=O(h^{-1}), C_{22}=O(1)$. All the numerical convergence rates of the exceptive cases are half-order higher than the theoretical convergence rates.
\begin{table}[htbp]
\centering
\tabcolsep 2pt
\caption{Unstructured triangular meshes: Error $\|\boldsymbol{u}-\boldsymbol{u}_h\|_0$ vs $h$
for different choices of $C_{11}, C_{22}$ when
$k=1, l=1$.}\label{table:errorL2k1l1unstructure} \vspace*{3pt}
\resizebox{0.98\textwidth}{!}{
\def\temptablewidth{\textwidth}
\begin{tabular*}{\temptablewidth}{@{\extracolsep{\fill}}|c|c|c|c|c|c|c|c|c|c|c|}
\hline \backslashbox{\kern3em$h$\kern-3em}{\kern2em$C_{11}, C_{22}$} & $O(h^{-1}), 0$ & order & $O(h^{-1}), O(1)$ & order & $O(h^{-1}), O(h)$ & order & $O(1), O(1)$ & order & $O(1), O(h)$ & order \\
\hline $1$ & 3.3990E+00 & $-$ & 3.9871E+00 & $-$ & 4.2302E+00 & $-$ & 3.8239E+00 & $-$ & 4.0069E+00 & $-$ \\
\hline $2^{-1}$ & 1.2308E+00 & 1.47 & 1.2207E+00 & 1.71 & 1.1359E+00 & 1.90 & 1.3569E+00 & 1.49 & 1.2274E+00 & 1.71 \\
\hline $2^{-2}$ & 4.5748E$-$01 & 1.43 & 3.1062E$-$01 & 1.97 & 2.8584E$-$01 & 1.99 & 3.6396E$-$01 & 1.90 & 3.0055E$-$01 & 2.03 \\
\hline $2^{-3}$ & 1.4422E$-$01 & 1.67 & 7.6164E$-$02 & 2.03 & 7.1621E$-$02 & 2.00 & 9.0756E$-$02 & 2.00 & 7.2383E$-$02 & 2.05 \\
\hline $2^{-4}$ & 4.0165E$-$02 & 1.84 & 1.8799E$-$02 & 2.02 & 1.7902E$-$02 & 2.00 & 2.2686E$-$02 & 2.00 & 1.7706E$-$02 & 2.03 \\
\hline $2^{-5}$ & 1.0503E$-$02 & 1.94 & 4.6480E$-$03 & 2.02 & 4.4728E$-$03 & 2.00 & 5.7187E$-$03 & 1.99 & 4.3599E$-$03 & 2.02 \\
\hline
       \end{tabular*}
       }
\end{table}
\begin{table}[htbp]
\centering
\tabcolsep 2pt \caption{Unstructured triangular meshes: Error $\left|(\boldsymbol{\sigma}-\boldsymbol{\sigma}_h,\boldsymbol{u}-\boldsymbol{u}_h)\right|_A$  vs $h$
for different choices of $C_{11}, C_{22}$ when
$k=1, l=1$.}\label{table:errorEnergyk1l1unstructure} \vspace*{3pt}
\resizebox{0.98\textwidth}{!}{
\def\temptablewidth{\textwidth}
\begin{tabular*}{\temptablewidth}{@{\extracolsep{\fill}}|c|c|c|c|c|c|c|c|c|c|c|}
\hline \backslashbox{\kern3em$h$\kern-3em}{\kern2em$C_{11}, C_{22}$} & $O(h^{-1}), 0$ & order & $O(h^{-1}), O(1)$ & order & $O(h^{-1}), O(h)$ & order & $O(1), O(1)$ & order & $O(1), O(h)$ & order \\
\hline $1$ & 9.2989E+00 & $-$ & 1.1613E+01 & $-$ & 1.1982E+01 & $-$ & 1.1572E+01 & $-$ & 1.1896E+01 & $-$ \\
\hline $2^{-1}$ & 5.6142E+00 & 0.73 & 5.6688E+00 & 1.03 & 5.5647E+00 & 1.11 & 5.7608E+00 & 1.01 & 5.5708E+00 & 1.09 \\
\hline $2^{-2}$ & 3.4557E+00 & 0.70 & 2.0151E+00 & 1.49 & 2.2948E+00 & 1.28 & 2.0813E+00 & 1.47 & 2.0977E+00 & 1.41 \\
\hline $2^{-3}$ & 1.9561E+00 & 0.82 & 6.6768E$-$01 & 1.59 & 1.0358E+00 & 1.15 & 7.2181E$-$01 & 1.53 & 8.5012E$-$01 & 1.30 \\
\hline $2^{-4}$ & 1.0372E+00 & 0.92 & 2.2314E$-$01 & 1.58 & 4.9826E$-$01 & 1.06 & 2.5343E$-$01 & 1.51 & 3.7363E$-$01 & 1.19 \\
\hline $2^{-5}$ & 5.3245E$-$01 & 0.96 & 7.5818E$-$02 & 1.56 & 2.4582E$-$01 & 1.02 & 8.9960E$-$02 & 1.49 & 1.7292E$-$01 & 1.11 \\
\hline
       \end{tabular*}
       }
\end{table}

We next consider the uniform rectangular meshes.
The numerical errors $\|\boldsymbol{u}-\boldsymbol{u}_h\|_0$ and $\left|(\boldsymbol{\sigma}-\boldsymbol{\sigma}_h,\boldsymbol{u}-\boldsymbol{u}_h)\right|_A$ for $k=1, l=1$ are listed in Tables~\ref{table:errorL2k1l1rectangle}-\ref{table:errorEnergyk1l1rectangle}, respectively.
We can see that $\|\boldsymbol{u}-\boldsymbol{u}_h\|_0=O(h^2)$ for different choices of $C_{11}$ and $C_{22}$, and
$\left|(\boldsymbol{\sigma}-\boldsymbol{\sigma}_h,\boldsymbol{u}-\boldsymbol{u}_h)\right|_A=O(h^{3/2})$ for $C_{11}=O(1), C_{22}=O(1)$ and
$\left|(\boldsymbol{\sigma}-\boldsymbol{\sigma}_h,\boldsymbol{u}-\boldsymbol{u}_h)\right|_A=O(h)$ for other choices of $C_{11}$ and $C_{22}$.
Therefore,
the numerical convergence rates of $\|\boldsymbol{u}-\boldsymbol{u}_h\|_0$  coincide with the theoretical results except the cases $C_{11}=O(h^{-1}), C_{22}=O(1)$ and $C_{11}=O(1), C_{22}=O(h)$. The numerical convergence rates of $\left|(\boldsymbol{\sigma}-\boldsymbol{\sigma}_h,\boldsymbol{u}-\boldsymbol{u}_h)\right|_A$ coincide with the theoretical results for different choices of $C_{11}$ and $C_{22}$.
\begin{table}[htbp]
\centering
\tabcolsep 2pt
\caption{Uniform rectangular meshes: Error $\|\boldsymbol{u}-\boldsymbol{u}_h\|_0$ vs $h$
for different choices of $C_{11}, C_{22}$ when
$k=1, l=1$.}\label{table:errorL2k1l1rectangle} \vspace*{3pt}
\resizebox{0.98\textwidth}{!}{
\def\temptablewidth{\textwidth}
\begin{tabular*}{\temptablewidth}{@{\extracolsep{\fill}}|c|c|c|c|c|c|c|c|c|c|c|}
\hline \backslashbox{\kern3em$h$\kern-3em}{\kern2em$C_{11}, C_{22}$} & $O(h^{-1}), 0$ & order & $O(h^{-1}), O(1)$ & order & $O(h^{-1}), O(h)$ & order & $O(1), O(1)$ & order & $O(1), O(h)$ & order \\
\hline $1$ & 3.7211E+00 & $-$ & 5.1403E+00 & $-$ & 5.1403E+00 & $-$ & 5.1403E+00 & $-$ & 5.1403E+00 & $-$ \\
\hline $2^{-1}$ & 1.4032E+00 & 1.41 & 1.2339E+00 & 2.06 & 1.2563E+00 & 2.03 & 1.5629E+00 & 1.72 & 1.5157E+00 & 1.76 \\
\hline $2^{-2}$ & 5.0795E$-$01 & 1.47 & 3.3884E$-$01 & 1.86 & 3.6796E$-$01 & 1.77 & 4.0140E$-$01 & 1.96 & 3.9904E$-$01 & 1.93 \\
\hline $2^{-3}$ & 1.6044E$-$01 & 1.66 & 9.1893E$-$02 & 1.88 & 1.0227E$-$01 & 1.85 & 9.7250E$-$02 & 2.05 & 1.0385E$-$01 & 1.94 \\
\hline $2^{-4}$ & 4.4923E$-$02 & 1.84 & 2.3921E$-$02 & 1.94 & 2.6772E$-$02 & 1.93 & 2.3343E$-$02 & 2.06 & 2.7241E$-$02 & 1.93 \\
\hline
       \end{tabular*}
       }
\end{table}
\begin{table}[htbp]
\centering
\tabcolsep 2pt \caption{Uniform rectangular meshes: Error $\left|(\boldsymbol{\sigma}-\boldsymbol{\sigma}_h,\boldsymbol{u}-\boldsymbol{u}_h)\right|_A$  vs $h$
for different choices of $C_{11}, C_{22}$ when
$k=1, l=1$.}\label{table:errorEnergyk1l1rectangle} \vspace*{3pt}
\resizebox{0.98\textwidth}{!}{
\def\temptablewidth{\textwidth}
\begin{tabular*}{\temptablewidth}{@{\extracolsep{\fill}}|c|c|c|c|c|c|c|c|c|c|c|}
\hline \backslashbox{\kern3em$h$\kern-3em}{\kern2em$C_{11}, C_{22}$} & $O(h^{-1}), 0$ & order & $O(h^{-1}), O(1)$ & order & $O(h^{-1}), O(h)$ & order & $O(1), O(1)$ & order & $O(1), O(h)$ & order \\
\hline $1$ & 1.2537E+01 & $-$ & 1.6321E+01 & $-$ & 1.6321E+01 & $-$ & 1.6321E+01 & $-$ & 1.6321E+01 & $-$ \\
\hline $2^{-1}$ & 6.4317E+00 & 0.96 & 6.9735E+00 & 1.23 & 6.8336E+00 & 1.26 & 7.0660E+00 & 1.21 & 6.8139E+00 & 1.26 \\
\hline $2^{-2}$ & 3.8075E+00 & 0.76 & 3.4415E+00 & 1.02 & 3.3483E+00 & 1.03 & 2.9002E+00 & 1.28 & 2.6527E+00 & 1.36 \\
\hline $2^{-3}$ & 2.2104E+00 & 0.79 & 1.6759E+00 & 1.04 & 1.6758E+00 & 1.00 & 1.1220E+00 & 1.37 & 1.0334E+00 & 1.36 \\
\hline $2^{-4}$ & 1.2003E+00 & 0.88 & 8.0902E$-$01 & 1.05 & 8.3997E$-$01 & 1.00 & 4.2251E$-$01 & 1.41 & 4.3790E$-$01 & 1.24 \\
\hline
       \end{tabular*}
       }
\end{table}

Finally let us verify the $p$-version convergence of our DG method on an uniform triangular mesh with fixed $h$ and take $\alpha_1=\beta_1=0$. We choose $h=1/8$ in this part, and take $l=k$. The errors $\|\boldsymbol{u}-\boldsymbol{u}_h\|_0$ and $\left|(\boldsymbol{\sigma}-\boldsymbol{\sigma}_h,\boldsymbol{u}-\boldsymbol{u}_h)\right|_A$ for $C_{11}=O(p)$, $C_{22}=O(1)$, $C_{11}=O(p)$, $C_{22}=O(1/p)$, $C_{11}=O(1)$, $C_{22}=O(1)$, and $C_{11}=O(1)$, $C_{22}=O(1/p)$ are listed in 
Tables~\ref{table:errorh18C11m1C220}-\ref{table:errorh18C110C22p1},
respectively.
As the $h$-version convergence numerically, we also observe that $\|\boldsymbol{u}-\boldsymbol{u}_h\|_0=O(1/p^{k+1})$ for different choices of $C_{11}$ and $C_{22}$, and $\left|(\boldsymbol{\sigma}-\boldsymbol{\sigma}_h,\boldsymbol{u}-\boldsymbol{u}_h)\right|_A=O(1/p^k)$ for $C_{22}=O(1/p)$ and $\left|(\boldsymbol{\sigma}-\boldsymbol{\sigma}_h,\boldsymbol{u}-\boldsymbol{u}_h)\right|_A=O(1/p^{k+0.5})$ for $C_{22}=O(1)$. Thus,
the numerical convergence rates of $\|\boldsymbol{u}-\boldsymbol{u}_h\|_0$  coincide with the theoretical results in Theorem~\ref{thm:5} except the cases $C_{11}=O(p), C_{22}=O(1)$ and $C_{11}=O(1), C_{22}=O(1/p)$, in which the numerical convergence rates of $\|\boldsymbol{u}-\boldsymbol{u}_h\|_0$ are half-order higher than the theoretical convergence rates. The numerical convergence rates of $\left|(\boldsymbol{\sigma}-\boldsymbol{\sigma}_h,\boldsymbol{u}-\boldsymbol{u}_h)\right|_A$  also coincide with the theoretical results in Theorem~\ref{thm:4} except the case $C_{11}=O(p), C_{22}=O(1)$, in which the numerical convergence rate of $\left|(\boldsymbol{\sigma}-\boldsymbol{\sigma}_h,\boldsymbol{u}-\boldsymbol{u}_h)\right|_A$ is half-order higher than the theoretical convergence rate.
\begin{table}[htbp]
\centering
\tabcolsep 2pt \caption{Uniform triangular meshes: Errors  vs $1/p$
when $C_{11}=O(p), C_{22}=O(1)$ and $l=k$.}\label{table:errorh18C11m1C220} \vspace*{3pt}
\def\temptablewidth{8cm}
\begin{tabular*}{\temptablewidth}{@{\extracolsep{\fill}}|c|c|c|}
\hline \quad $k$\;\;\;  & $p^{k+1}\|\boldsymbol{u}-\boldsymbol{u}_h\|_0$ & $p^{k+0.5}\left|(\boldsymbol{\sigma}-\boldsymbol{\sigma}_h,\boldsymbol{u}-\boldsymbol{u}_h)\right|_A$ \\
\hline \quad 1\;\;\;  & 0.2812 & 1.6906 \\
\hline \quad 2\;\;\;  & 0.0704 & 0.5717 \\
\hline \quad 3\;\;\;  & 0.0227 & 0.1782 \\
\hline \quad 4\;\;\;  & 0.0085 & 0.0730 \\
\hline \quad 5\;\;\;  & 0.0107 & 0.1352 \\
\hline
       \end{tabular*}
\end{table}
\begin{table}[htbp]
\centering
\tabcolsep 2pt \caption{Uniform triangular meshes: Errors  vs $1/p$
when $C_{11}=O(p), C_{22}=O(1/p)$ and $l=k$.}\label{table:errorh18C11m1C22p1} \vspace*{3pt}
\def\temptablewidth{7.5cm}
\begin{tabular*}{\temptablewidth}{@{\extracolsep{\fill}}|c|c|c|}
\hline \quad $k$\;\;\;  & $p^{k+1}\|\boldsymbol{u}-\boldsymbol{u}_h\|_0$ & $p^{k}\left|(\boldsymbol{\sigma}-\boldsymbol{\sigma}_h,\boldsymbol{u}-\boldsymbol{u}_h)\right|_A$ \\
\hline \quad 1\;\;\;  & 0.2812 & 1.1954 \\
\hline \quad 2\;\;\;  & 0.0620 & 0.3244 \\
\hline \quad 3\;\;\;  & 0.0202 & 0.0984 \\
\hline \quad 4\;\;\;  & 0.0077 & 0.0376 \\
\hline \quad 5\;\;\;  & 0.0086 & 0.0522 \\
\hline
       \end{tabular*}
\end{table}
\begin{table}[htbp]
\centering
\tabcolsep 2pt \caption{Uniform triangular meshes: Errors  vs $1/p$
when $C_{11}=O(1), C_{22}=O(1)$ and $l=k$.}\label{table:errorh18C110C220} \vspace*{3pt}
\def\temptablewidth{8cm}
\begin{tabular*}{\temptablewidth}{@{\extracolsep{\fill}}|c|c|c|}
\hline \quad $k$\;\;\;  & $p^{k+1}\|\boldsymbol{u}-\boldsymbol{u}_h\|_0$ & $p^{k+0.5}\left|(\boldsymbol{\sigma}-\boldsymbol{\sigma}_h,\boldsymbol{u}-\boldsymbol{u}_h)\right|_A$ \\
\hline \quad 1\;\;\;  & 0.2812 & 1.6906 \\
\hline \quad 2\;\;\;  & 0.0794 & 0.5119 \\
\hline \quad 3\;\;\;  & 0.0260 & 0.1594 \\
\hline \quad 4\;\;\;  & 0.0092 & 0.0532 \\
\hline \quad 5\;\;\;  & 0.0125 & 0.1449 \\
\hline
       \end{tabular*}
\end{table}
\begin{table}[htbp]
\centering
\tabcolsep 2pt \caption{Uniform triangular meshes: Errors  vs $1/p$
when $C_{11}=O(1), C_{22}=O(1/p)$ and $l=k$.}\label{table:errorh18C110C22p1} \vspace*{3pt}
\def\temptablewidth{7.5cm}
\begin{tabular*}{\temptablewidth}{@{\extracolsep{\fill}}|c|c|c|}
\hline \quad $k$\;\;\;  & $p^{k+1}\|\boldsymbol{u}-\boldsymbol{u}_h\|_0$ & $p^{k}\left|(\boldsymbol{\sigma}-\boldsymbol{\sigma}_h,\boldsymbol{u}-\boldsymbol{u}_h)\right|_A$ \\
\hline \quad 1\;\;\;  & 0.2812 & 1.1954 \\
\hline \quad 2\;\;\;  & 0.0659 & 0.2800 \\
\hline \quad 3\;\;\;  & 0.0210 & 0.0811 \\
\hline \quad 4\;\;\;  & 0.0077 & 0.0266 \\
\hline \quad 5\;\;\;  & 0.0098 & 0.0552 \\
\hline
       \end{tabular*}
\end{table}

\subsection{A three-dimensional example}

Let
$\Omega$ be the unit cube $(0,1)^3$,  and
\begin{align*}
\boldsymbol{u}(x_1, x_2, x_3)=&\left(\begin{array}{l}
\sin(\pi x_1)\sin(\pi x_2)\sin(\pi x_3) \\
2\sin(\pi x_1)\sin(\pi x_2)\sin(\pi x_3) \\
4\sin(\pi x_1)\sin(\pi x_2)\sin(\pi x_3)
\end{array}\right).
\end{align*}
The right hand side $\boldsymbol{f}$ is computed from problem \eqref{eq:elas}.
We adopt the uniform tetrahedral meshes $\mathcal{T}_h$ of $\Omega$.
The numerical errors $\|\boldsymbol{u}-\boldsymbol{u}_h\|_0$ and $\left|(\boldsymbol{\sigma}-\boldsymbol{\sigma}_h,\boldsymbol{u}-\boldsymbol{u}_h)\right|_A$ for $k=1, l=1$ are shown in
Tables~\ref{table:errorL2k1l13d}-\ref{table:errorEnergyk1l13d}, respectively,
from which we observe that $\|\boldsymbol{u}-\boldsymbol{u}_h\|_0=O(h^2)$ for different choices of $C_{11}$ and $C_{22}$, and $\left|(\boldsymbol{\sigma}-\boldsymbol{\sigma}_h,\boldsymbol{u}-\boldsymbol{u}_h)\right|_A=O(h)$ for $C_{22}=O(h), 0$ and $\left|(\boldsymbol{\sigma}-\boldsymbol{\sigma}_h,\boldsymbol{u}-\boldsymbol{u}_h)\right|_A=O(h^{3/2})$ for $C_{22}=O(1)$. Hence the numerical convergence rates of $\|\boldsymbol{u}-\boldsymbol{u}_h\|_0$ and $\left|(\boldsymbol{\sigma}-\boldsymbol{\sigma}_h,\boldsymbol{u}-\boldsymbol{u}_h)\right|_A$ are same as those of the two-dimensional example.
To be specific,
the numerical convergence rates of $\|\boldsymbol{u}-\boldsymbol{u}_h\|_0$ coincide with the theoretical results except the cases $C_{11}=O(h^{-1}), C_{22}=O(1)$ and $C_{11}=O(1), C_{22}=O(h)$. And the numerical convergence rates of $\left|(\boldsymbol{\sigma}-\boldsymbol{\sigma}_h,\boldsymbol{u}-\boldsymbol{u}_h)\right|_A$ coincide with the theoretical results except the case $C_{11}=O(h^{-1}), C_{22}=O(1)$.
\begin{table}[htbp]
\centering
\tabcolsep 2pt
\caption{Uniform tetrahedral meshes: Error $\|\boldsymbol{u}-\boldsymbol{u}_h\|_0$ vs $h$
for different choices of $C_{11}, C_{22}$ when
$k=1, l=1$.}\label{table:errorL2k1l13d} \vspace*{3pt}
\resizebox{0.98\textwidth}{!}{
\def\temptablewidth{\textwidth}
\begin{tabular*}{\temptablewidth}{@{\extracolsep{\fill}}|c|c|c|c|c|c|c|c|c|c|c|}
\hline \backslashbox{\kern3em$h$\kern-3em}{\kern2em$C_{11}, C_{22}$} & $O(h^{-1}), 0$ & order & $O(h^{-1}), O(1)$ & order & $O(h^{-1}), O(h)$ & order & $O(1), O(1)$ & order & $O(1), O(h)$ & order \\
\hline $2^{-1}$ & 4.1570E$-$01 & $-$ & 4.2975E$-$01 & $-$ & 3.7248E$-$01 & $-$ & 5.5379E$-$01 & $-$ & 3.9309E$-$01 & $-$ \\
\hline $2^{-2}$ & 1.3978E$-$01 & 1.57 & 1.3181E$-$01 & 1.71 & 1.0696E$-$01 & 1.80 & 1.5881E$-$01 & 1.80 & 9.3257E$-$02 & 2.08 \\
\hline $2^{-3}$ & 4.0197E$-$02 & 1.80 & 3.1358E$-$02 & 2.07 & 2.7547E$-$02 & 1.96 & 3.8824E$-$02 & 2.03 & 2.2015E$-$02 & 2.08 \\
\hline $2^{-4}$ & 1.0654E$-$02 & 1.92 & 7.5727E$-$03 & 2.05 & 6.9252E$-$03 & 1.99 & 9.4496E$-$03 & 2.04 & 5.3621E$-$03 & 2.04 \\
\hline
       \end{tabular*}
       }
\end{table}
\clearpage \begin{table}[htbp]
\centering
\tabcolsep 2pt \caption{Uniform tetrahedral meshes: Error $\left|(\boldsymbol{\sigma}-\boldsymbol{\sigma}_h,\boldsymbol{u}-\boldsymbol{u}_h)\right|_A$  vs $h$
for different choices of $C_{11}, C_{22}$ when
$k=1, l=1$.}\label{table:errorEnergyk1l13d} \vspace*{3pt}
\resizebox{0.98\textwidth}{!}{
\def\temptablewidth{\textwidth}
\begin{tabular*}{\temptablewidth}{@{\extracolsep{\fill}}|c|c|c|c|c|c|c|c|c|c|c|}
\hline \backslashbox{\kern3em$h$\kern-3em}{\kern2em$C_{11}, C_{22}$} & $O(h^{-1}), 0$ & order & $O(h^{-1}), O(1)$ & order & $O(h^{-1}), O(h)$ & order & $O(1), O(1)$ & order & $O(1), O(h)$ & order \\
\hline $2^{-1}$ & 3.1104E+00 & $-$ & 3.3806E+00 & $-$ & 3.1350E+00 & $-$ & 3.5101E+00 & $-$ & 3.0264E+00 & $-$ \\
\hline $2^{-2}$ & 1.7711E+00 & 0.81 & 1.2862E+00 & 1.39 & 1.3598E+00 & 1.21 & 1.3165E+00 & 1.41 & 1.1027E+00 & 1.46 \\
\hline $2^{-3}$ & 9.5073E$-$01 & 0.90 & 4.4060E$-$01 & 1.55 & 6.2062E$-$01 & 1.13 & 4.5439E$-$01 & 1.53 & 4.2521E$-$01 & 1.37 \\
\hline $2^{-4}$ & 4.8999E$-$01 & 0.96 & 1.6657E$-$01 & 1.40 & 3.0027E$-$01 & 1.05 & 1.5741E$-$01 & 1.53 & 1.8309E$-$01 & 1.22 \\
\hline
       \end{tabular*}
       }
\end{table}

\section{Conclusion}

The $hp$-version error analysis is systematically developed for the mixed DG method  \eqref{eq:variation} (or equivalently, the method \eqref{eq:mix1}-\eqref{eq:mix2}). The derivation is mainly based on the ideas in \cite{CastilloCockburnPerugiaSchotzau2000} and the $hp$-version error estimates of two $L^2$ projection operators. According to our numerical experiments, we may achieve the following conclusions:
\begin{enumerate}[(1)]
\item If $k=l+1$, the error estimates of $\|\boldsymbol{u}-\boldsymbol{u}_h\|_0$ and $\left|(\boldsymbol{\sigma}-\boldsymbol{\sigma}_h,\boldsymbol{u}-\boldsymbol{u}_h)\right|_A$ in Theorems~\ref{thm:4}-\ref{thm:5} are sharp except the case $(C_{11}, C_{22})=(O(1), O(h))$.
\item If $k=l>0$, the error estimate of $\left|(\boldsymbol{\sigma}-\boldsymbol{\sigma}_h,\boldsymbol{u}-\boldsymbol{u}_h)\right|_A$ in Theorem~\ref{thm:4} is sharp for all choices of $C_{11}$ and $C_{22}$.
\item If $k=l>0$, the error estimate of $\|\boldsymbol{u}-\boldsymbol{u}_h\|_0$ in Theorem~\ref{thm:5} is sharp except the cases $(C_{11}, C_{22})=(O(h^{-1}), O(1))$ and $(C_{11}, C_{22})=(O(1), O(h))$.
\item All the numerical convergence rates are half-order higher than the theoretical convergence rates if the theoretical error estimates are not sharp.
\end{enumerate}
According to the error analysis developed in the last section, we can not show the last point of the conclusions in theory. It is our further work to study the sharp theoretical error estimate of $\|\boldsymbol{u}-\boldsymbol{u}_h\|_0$ for the cases $(C_{11}, C_{22})=(O(h^{-1}), O(1))$ and $(C_{11}, C_{22})=(O(1), O(h))$.

The other interesting and valuable work is to extend the mixed DG method \eqref{eq:variation} for numerically solving linear transient elasticity problems as well as nonlinear elasticity problems. And then develop the $hp$-version error analysis for the corresponding numerical methods.



\end{document}